\documentclass{article}
\usepackage{amsmath,amssymb,amscd,stmaryrd}
\usepackage[all]{xy}
\usepackage{caption,graphicx}

\font\ibf=cmbxti10

\title{On the 3-Dimensional Poincar\'e Conjecture and the 4-Dimensional Smooth Schoenflies Problem \\
{\Large -- a double research announcement --}}
\author{\it by \\
{ \ } \\
\Large Valentin Po\'enaru\footnote{Universit\'e de Paris-Sud, Math\'ematiques 425, Topologie et Dynamique, 91405 Orsay Cedex, France. 
\newline
Universit\'a degli Studi di Trento, Dipartimento di Matematica, 38050 Povo-Trento, Italia. 
\newline
This paper has been partially supported by the NSF grant DMS-0071852.}}
\date{July 2006}

\begin{document}

\thispagestyle{empty}

\maketitle

\newpage

\noindent {\Large \bf 0 \quad A foreword}

\medskip

\noindent This is the announcement of an alternative approach to the 3-dimensional Poincar\'e Conjecture, different from Perelman's big and spectacular breakthrough. No claim concerning the other parts of the Thurston Geometrization Conjecture, come with our purely 4-dimensional line of argument.

\smallskip

The format of the present paper is that of a rather informal letter which, conceivably, I might have written to some mathematical friend, let us say to David Gabai.

\vglue 1cm

\section{Introduction}\label{sec1}

The present shortish paper is a double research announcement. On the one hand, my 4-dimensional program for proving the 3-dimensional Poincar\'e Conjecture, developed during several decades, is now finally completely finished. The last step which was missing is provided by the Theorem~1 below. But then also, there is an outgrowth of this program, namely the proof that any smooth 4-dimensional Schoenflies ball is geometrically simply-connected, i.e. it possesses smooth handlebody decompositions without handles of index one. This is Theorem~2 below.

\smallskip

My almost purely 4-dimensional techniques for the Poincar\'e Conjecture, are completely independent, of course, of the Ricci flow approach of R.~Hamilton and G.~Perelman; see here \cite{Mo} for more extensive references.

\smallskip

It is the two theorems 1 and 2 mentioned above and then stated precisely below, which are the novelties here, with respect to my 2004 informal and more tentative announcement in the Steklov Proceedings \cite{Po-S}.

\bigskip

\noindent {\bf Theorem 1.} -- {\it For any homotopy $3$-ball $\Delta^3$, we introduce the following, canonically attached, open smooth $4$-manifold
\begin{equation}
\label{eq1.1}
X^4 = {\rm int} \, [(\Delta^3 \times I) \, \# \, \infty \, \# \, (S^2 \times D^2)] \, .
\end{equation}
IF this $X^4$ is geometrically simply connected, THEN so is $\Delta^3 \times I$ itself.}

\bigskip

So, the format of this statement is, actually, the following implication
\begin{equation}
\label{eq1.2}
X^4 ({\rm open}) \ {\rm g.s.c.} \Longrightarrow \Delta^3 \times I ({\rm compact}) \ {\rm g.s.c.}
\end{equation}

Next, we consider Schoenflies 4-balls. By definition these are smooth compact 4-manifolds, which we will denote generically by $\Delta_{\rm Schoenflies}^4$, such that
$$
\partial \Delta_{\rm Schoenflies}^4 = S^3, \ \mbox{and there is a smooth embedding} \ \Delta_{\rm Schoenflies}^4 \subset S_{\rm standard}^4 \, .
$$
With this, here is our

\bigskip

\noindent {\bf Theorem 2.} -- {\it Any $\Delta_{\rm Schoenflies}^4$ is geometrically simply connected.}

\bigskip

Remember here that, according to the classical work of Barry Mazur \cite{Ma}, it is certainly known that such a $\Delta_{\rm Sch}^4$ is topologically the 4-ball. Even better, if we delete from it a boundary point, then what we get is diffeomorphic to the standard 4-ball, with a boundary point removed.

\smallskip

The plan of the present paper is the following. In the next section 2 we will give a bird's eye view short outline of my proof of the Poincar\'e Conjecture, showing in particular how Theorem 1 above fits into it. A much more detailed outline, but with Theorem 1 occuring there with a question mark, was given in the Steklov paper \cite{Po-S}, and so we will be really very brief here. In the same \cite{Po-S}, which may be considered a companion of the present announcement, Theorem 2 was only hinted at, as a possibility.

\smallskip

Next, in the sections 3 and 4, we will give a glimpse of how the proofs of Theorems 1 and 2 go. As it stands, section 3 should give already a very first idea, while the even more impressionistic section 4 touches on some more technical issues.

\smallskip

But the point is that both proofs can and will be presented, largely, simultaneously, and in the same breath. Here is how one should view their starting points. For Theorem 1, the starting point is an {\it hypothesis}, namely the $X^4$ (\ref{eq1.1}) being g.s.c.; Theorem 2 starts from the fact that, as a consequence of Barry's work mentioned above, the open smooth 4-manifold $\Delta^4 \cup (S^3 \times [0,\infty))$ ($= {\rm int} \, \Delta^4$), which we may as well call now $X^4$ again, {\ibf is} geometrically simply connected. Of course, Barry's work really implies that it is the standard $R^4$. But only the g.s.c. property will be retained here, for our present purposes. I do believe that, afterwards, in order to show that any g.s.c. $\Delta_{\rm Schoenflies}^4$ is actually standard, the full strength of Barry's result (as far as the $4^d$ DIFF category goes) should be needed. But this is another story.

\smallskip

The complete detailed proofs of Theorems 1 and 2 are completely (hand-) written down, in a very long two-part paper, to which I will refer hereafter as ``PoV-B'', listed as \cite{PoV-B}. I hope to be able to make it available in a typed version in a not too long time. I should also add that working on Theorem~2 has been, for me, a very good testing ground for the PoV-B technology of Theorem~1, with a lot of feed-back between the two items.

\smallskip

I owe too much to too many people to start listing them all here and now. This notwithstanding, I do want to thank David Gabai, without the help of whom this work would not have been here. And then, I should also mention that the very first impetus for trying to link together things like in the two theorems above, came from a suggestion which Michael Freedman has made to David and me, way back in the Spring 1995.

\smallskip

Finally, thanks are due to the IH\'ES for generously offering me the possibility to use its typing facilities and to C\'ecile Cheikhchoukh and Marie-Claude Vergne for the typing and the drawings.

\section{A brief outline of the proof of the Poincar\'e Conjecture}\label{sec2}
\setcounter{equation}{0}

There are three distinct steps in the proof, each short to state but also each with a very long proof. I will present them here as follows.

\bigskip

\noindent STEP I. Here, the climax is the following final result

\bigskip

\noindent {\bf Theorem 3.} -- {\it For any homotopy $3$-ball $\Delta^3$, the open smooth $4$-manifold $X^4$ from {\rm (\ref{eq1.1})} is geometrically simply connected.}

\bigskip

The complete detailed proof is contained in the series of papers \cite{PoI}, \cite{PoII}, \cite{PoIII}, \cite{PoIV-A}, \cite{PoIV-B} and \cite{PoV-A}. For this last paper, which proves a 4-dimensional result, completely independent of the rest, there is also a shorter version \cite{PoV-Aoutline}. For this and the next step, see also \cite{Ga}, \cite{Po-B}, \cite{Po-T}. Notice that what Theorem 3 does, is to prove the hypothesis occuring in Theorem 1, i.e. in the implication (\ref{eq1.2}).

\bigskip

\noindent STEP II. Just like it was already the case for Theorem 1, the main result takes here again the form of an implication, namely, now the following
\begin{equation}
\label{eq2.1}
\Delta^3 \times I \ \mbox{geometrically simply connected} \ \Longrightarrow \Delta^3 \ {\rm standard}.
\end{equation}
More explicitly, this is the following

\bigskip

\noindent {\bf Theorem 4.} -- {\it Let $\Delta^3$ be a homotopy $3$-sphere which is such that $\Delta^3 \times I$ is geometrically simply connected. Then $\Delta^3 = B^3$.}

\bigskip

A brief outline of the proof can be found in \cite{Po-S}, \cite{Po-B}. The complete detailed proof is contained in \cite{PoVI}. 

\smallskip

Here are some words concerning the way in which the hypothesis that 
$\Delta^3 \times I$ is g.s.c. is being used in the proof of Theorem 4. At this point here is a little lemma; the various terms which are used in the statement are all explained in \cite{Ga} or \cite{PoI} (and \cite{PoII}).

\bigskip

\noindent {\bf Lemma 5.} -- {\it Let $\Delta^3$ be a homotopy $3$-ball which is such that $\Delta^3 \times I$ is g.s.c. (in the smooth category). Then there exists a collapsible pseudo-spine representation for $\Delta^3$, call it $K^2 \overset{f}{\longrightarrow} \Delta^3$, for which one can find a desingularization $\varphi$ having the following property. There exists a strategy for zipping $f$, which is COHERENT for $\varphi$.}

\bigskip

What ``coherence'' means here is that when, during the zipping process any two singularities $s_1 , s_2$ meet in a head-on collision, then their desingularizations are well-matched together: the $S(N)$ branches of $s_1$ match the $S(N)$ branches of $s_2$ (and not the $N(S)$ branches). The proof of Lemma 5 can be found in \cite{Ga}, \cite{PoV-A} and the converse statement to Lemma 5, going from coherence to g.s.c. is true too; actually it is even easier (see \cite{Ga} and \cite{PoII}).

\smallskip

Now, the point is that the starting point of the infinite processes via which Theorem 4 is proved, is a collapsible pseudo-spine representation for $\Delta^3$, having the coherence property. This is how the geometric simple connectivity of $\Delta^3 \times I$ comes in.

\bigskip

\noindent STEP III. This step is our present Theorem 1. 

\bigskip

The ordering of our three steps above was chronology rather than logic. This being said, on the same lines as the references \cite{PoI} to \cite{PoIV-B} and \cite{PoV-A}, all part of step I above, just after them and before \cite{PoVI} (step II), there should actually come now the \cite{PoV-B} of step III, containing the proof of Theorem 1 and, incidentally, of Theorem 2 too.

\smallskip

Once one assumes all the three steps above, one can plug Theorem 3 (step I) into Theorem 1 (step III) and conclude that $\Delta^3 \times I$ is always geometrically simply connected. When this fact is plugged into Theorem 4 (step II), then this yields the following main result

\bigskip

\noindent {\bf Theorem 6.} (THE POINCAR\'E CONJECTURE) -- {\it Every homotopy $3$-ball $\Delta^3$ is standard, i.e. $\Delta^3 = B^3$.}

\bigskip

Notice that, on the way, we have also proved the so-called COHERENCE THEOREM, stating that {\it every} homotopy 3-ball admits a collapsible pseudo-spine representations which is also coherent. The gap in an earlier attempted direct proof for the coherence theorem, which was detected in the Spring 1995 by Michael Freedman and David Gabai, is now completely filled in by the combination of \cite{PoV-A} (which in the presentation chosen here has been included inside Step I) together with the Theorem 1 above, i.e. by \cite{PoV-B}. Together, the PoV-A (\cite{PoV-A}) and PoV-B (\cite{PoV-B}) completely supersede the by now dead Orsay preprint 94-25 \cite{PoV}, from 1994. One may also put these things as follows. The \cite{PoIV-B} proves that $\Delta^3 \times I$ has the property of being geometrically simply connected {\ibf at long distance} (this is a notion weaker than g.s.c., for which I refer to \cite{Po-S}, \cite{Po-B}). Then what \cite{PoV-A} $+$ \cite{PoV-B} actually do for us, is to deduce the COHERENCE THEOREM from this last property.

\section{Some hints concerning the proof of the Theorems 1 and 2}\label{sec3}
\setcounter{equation}{0}

We will try, as much as possible, to present the two proofs in parallel; the fact that we may do this, should be seen as a distinctive feature of the present approach. Everything now is in the smooth category and we will denote by $\Delta^4$ a compact bounded 4-manifold which is either $\Delta^3 \times I$ (with $\Delta^3$ a homotopy 3-ball) or $\Delta_{\rm Schoenflies}^4$.

\smallskip

We start with the following sequence of nested spaces
\begin{equation}
\label{eq3.1}
\Delta^4 = \Delta_{\rm small}^4 \subset X_{\rm open}^4 \subset \Delta_{\rm large}^4 = \Delta_1^4 \, ,
\end{equation}
where $\Delta_{\rm small}^4$ and $\Delta_{\rm large}^4$ are two copies of the same $\Delta^4$, separated by a product collar. The $X^4$ is, according to the case, either the ${\rm int} \, [(\Delta^3 \times I) \, \# \, \infty \, \# \, (S^2 \times D^2)]$ from (\ref{eq1.1}) or, in the Schoenflies case it is $\Delta^4 \cup (\partial \Delta^4 \times [0,\infty))$. In both cases we have a splitting of the form $X^4 = X^3 \times R$, but there is of course no compact splitting like $\Delta^3 \times I$, in the Schoenflies case. From the very beginning we are presented here with two distinct features or structures, referred hereafter as ``RED'' and ``BLUE'' respectively. The RED feature of (\ref{eq3.1}) is a ``collapse'' $X^4 \to \Delta^4$; but this requires some qualifications. On the one hand, since $X^4$ is not compact, we should rather talk about an infinite dilatation process, going the other way around. But then also, more seriously, in the case of $\Delta^4 = \Delta^3 \times I$, the collapse (and we drop the quotation marks from now on), certainly has some {\ibf defects}, namely the infinitely many $\# \, (S^2 \times D^2)$ of the corresponding $X^4$. In our present brief outline we will chose to rather ignore them. Of course, in a more realistic discussion they will have to be dealt with; but see here also the remark which follows after (\ref{eq3.33.2}). But then also, they are absent in the Schoenflies context.

\smallskip

The BLUE feature of (\ref{eq3.1}) is that $X^4$, as such, is geometrically simply connected; here the compact $\Delta^4$ is altogether being ignored. In a combinatorial language, $X^4$ admits a smooth cell-decomposition with a 2-skeleton which is
$$
\mbox{(a collapsible infinite 2-complex)} + \mbox{(2-cells added)}.
$$
What we may hope to achieve by combining the two features above would be to construct, inside the collar $\Delta_1^4 - {\rm int} \, \Delta^4$ coming with (\ref{eq3.1}), a system of embedded exterior discs, in cancelling position with the 1-handles of $\Delta^4$. It is not very hard to show that this would imply that $\Delta^4$ is g.s.c., i.e. we would get both Theorems 1 and 2 this way. Presumably, the exterior discs should be gotten starting from the BLUE structure, and the connection with our $\Delta^4$, i.e. the cancelling property should be gotten by invoking the RED feature too. This sounds more like a vague pipe-dream, of course, but it may still serve as a vague guide-line for what will be following next.

\smallskip

But before we can really start off the ground, we will need to change the initial set-up (\ref{eq3.1}), in several successive stages.

\bigskip

\noindent STAGE I. Let us start by denoting $\Delta^2$ the 2-spine of $\Delta^3 = \Delta^3 \times 0 \subset \Delta^3 \times I$, {\it or} the 2-skeleton of $\Delta_{\rm Schoenflies}^4$ (so as to avoid having to deal with the 3-handles of $\Delta_{\rm Schoenflies}^4$), according to the case. The point is that all we need is to show that the $4^d$ regular neighbourhood $N^4 (\Delta^2)$ is g.s.c. It is on $\Delta^2$, rather than on $\Delta^4$, that we will focus from now on. We may even call $N^4 (\Delta^2) , \Delta^4$.

\smallskip

Now, a priori both the BLUE and the RED features, are each expressible in terms of two smooth cell-decompositions of $X^4$, the two being independent of each other. But then, making use of the smooth Hauptvermutung of J.H.C.~Whitehead \cite{W1} and also of some combinatorial arguments, into which we will not go here, one can produce a unique smooth cell-decomposition of $X^4$ which exhibits both the BLUE and the RED features. This will be done in terms of some combinatorial data to be explained now; this may be a bit lengthy, but it is unavoidable for our exposition.

\smallskip

Let $X^2$ be the 2-skeleton of $X^4$, and let also $\Gamma (\infty) \subset X^2$ be the 1-skeleton. The $\Delta^2 \subset X^2$ is a subcomplex, with its finite 1-skeleton $\Gamma (1) \subset \Gamma (\infty)$. Inside $\Gamma (\infty)$ live two independent set of points
\begin{equation}
\label{eq3.2}
R \ (\mbox{for red}) \subset \Gamma (\infty) \supset B \ (\mbox{for blue}) \, .
\end{equation}
Both $\Gamma (\infty) - R$ and $\Gamma (\infty) - B$ are trees, making that a given edge $e \subset \Gamma (\infty)$ contains at most one $R_i \in R$ and one $B_j \in B$. When $R_i \in e \ni B_j$, it will be assumed that $R_i = B_j \in R \cap B$. One should think here of the $R_i , B_j$ as being 1-handles or, more precisely 1-handle cocores, i.e. properly embedded 3-balls $(B^3 , \partial B^3) \subset (N^4 (\Gamma (\infty) , \partial N^4 (\Gamma (\infty))$. For further purposes, the following notations will be introduced too
\begin{equation}
\label{eq3.3}
R \, \cap \, \Gamma (1) = \{ R_1 , R_2 , \ldots , R_n \} \ \mbox{and} \ R - \{ R_1 , R_2 , \ldots , R_n \} = \{ h_1 , h_2 , h_3 , \ldots \} \, .
\end{equation}
One gets the $X^2$ and/or the $N^4 (X^2)$ by adding 2-cells and/or 2-handles along an infinite framed link
$$
\{ \mbox{link} \} \subset \partial N^4 (\Gamma (\infty)) \approx \Gamma (\infty) \, .
$$
The link comes with two independent disjoined partitions
\begin{eqnarray}
\label{eq3.4}
\{ \mbox{link} \} &= &\sum_1^n \Gamma_i + \sum_1^{\infty} C_j + \sum_1^{\infty} \gamma_k^0 \quad (\mbox{RED partition}) \\
&= &\sum_1^{\infty} \eta_{\ell} + \sum_1^{\infty} \gamma_m^1 \quad (\mbox{BLUE partition}) \, . \nonumber
\end{eqnarray}
For each element of a link we have an associated 2-cell and/or 2-handle, denoted in both cases by $D^2$ (curve). With this, (\ref{eq3.4}) leads to the following decompositions
\begin{equation}
\label{eq3.5}
\Delta^2 = \Gamma (1) \cup \sum_1^n D^2 (\Gamma_i) \, , \quad \mbox{and}
\end{equation}
\begin{eqnarray}
X^2 &= &\Gamma (\infty) \cup \left( \sum_1^n D^2 (\Gamma_i) + \sum_1^{\infty} D^2 (C_j) + \sum_1^{\infty} D^2 (\gamma_k^0)\right) \nonumber \\
&= & \Gamma (\infty) \cup \left( \sum_1^{\infty} D^2 (\eta_{\ell}) + \sum_1^{\infty} D^2 (\gamma_m^1)\right) \, . \nonumber
\end{eqnarray}

\bigskip

\noindent {\bf Remark.} We have used the same $n$, which by all means will mean the cardinality of $R \cap \Gamma (1)$, for the cardinality of the 2-handles $D^2(\Gamma)$ of $\Delta^2$ too. Now, this is perfectly legitimate in the case $\Delta^3 \times I$, when $\Delta^2$ is the spine. In the Schoenflies case, $\Delta^2$ is the 2-skeleton and, then actually
$$
\bar n = {\rm card} \, (D^2 (\Gamma)) > n = {\rm card} \, (R \cap \Gamma (1)) \, .
$$
Once this is understood, there should be no problem concerning this ambiguity in notation.

\bigskip

So, each curve and each disc comes with two independent labels, a red one and a blue one. The $X^4$ comes with a big RED collapsing flow (with possible defects in the $\Delta^3 \times I$ case). With this, the $D^2 (\gamma_k^0)$ are essentially (but see here also the remark below) those 2-cells which are killed by the $3^d$ RED collapse, while the $D^2 (C_j)$ are the 2-cells killed by the $2^d$ RED collapse. [Remark. In the Schoenflies case, all the $D^2 (\gamma^0)$ are {\it rigorously} killed by the RED 3-dimensional flow. In the $\Delta^3 \times I$ case, the normal $D(\gamma^0)$'s are, but then we also have non-trivial $D^2 (\gamma^0)$'s corresponding to the defects. Quite some care has to be devoted to them in real life. But in this exposition we will largely ignore them (as much as that will be possible).]

\smallskip

In terms of (\ref{eq3.3}) and (\ref{eq3.4}) we define the red geometric intersection matrix $C \cdot h$. We express the RED 2-dimensional collapse by stipulating that $C \cdot h$ is of the following {\ibf easy id + nilpotent} form (after appropriate re-indexing)
\begin{equation}
\label{eq3.6}
C_i \cdot h_j = \delta_{ij} + \xi_{ij}^0 \, , \ \mbox{where we can have $\xi_{ij}^0 \ne 0$ only if $i > j$}Ê\, .
\end{equation}

Finally, the BLUE feature is expressed by stipulating that
\begin{equation}
\label{eq3.7}
\mbox{The blue geometric intersection matrix $\eta \cdot B$ is also of the easy id $+$ nilpotent form.} 
\end{equation}

With all this, in the framework of a common, unique smooth cell-decomposi\-tion for $X^4$, we have encoded in a convenient combinatorial language both our red and blue features.

\bigskip

\noindent {\bf Remark.} There is also a notion of {\ibf difficult} id $+$ nilpotent. With notations like in (\ref{eq3.6}), this means now that
$$
\xi_{ij}^0 \ne 0 \quad \mbox{only if} \quad i < j \, .
$$
This, contrary to the easy id $+$ nil, is very far from ``collapsible'', when we are in the infinite context. It can be shown without difficulty, that the classical Whitehead manifold ${\rm Wh}^3$ \cite{W2} admits a handlebody decomposition with only handles of index one and two and with a geometric intersection which is of this type. But then, in \cite{PoV-A} (and see \cite{PoV-Aoutline} too) where it occurs quite naturally, the difficult id $+$ nil turned out to be quite useful too.

\smallskip

In the set-up which we have just introduced, we have the following two very useful objects
\begin{equation}
\label{eq3.8}
X_0^2 \underset{\rm def}{=} \Gamma (\infty) \cup \left( \sum_1^n D^2 (\Gamma_i) + \sum_1^{\infty} D^2 (C_j)\right) \supset \Gamma (\infty) \cup \sum_{1}^{\infty} D^2 (C_j) \, .
\end{equation}
There is now a RED 2-dimensional collapse of $X_0^2$ onto $\Delta^2$. But bluewise, the $X_0^2$ is clearly limping. From now on, we take $N^4 (\Delta^2)$ as being $\Delta^4$ and, correspondingly, $N^4 (\Delta^2) \cup \{\mbox{collar}\}$ as being $\Delta_1^4$. With all this, we change the set-up (\ref{eq3.1}) into the following, for the time being
\begin{equation}
\label{eq3.9}
\Delta^4 \subset N^4 (X_0^2) \subset N^4 (X^2) \subset \Delta_1^4 \, .
\end{equation}

\bigskip

\noindent STAGE II. This will be, essentially, a refinement of the previous stage, in preparation for the next things to come. In the context of (\ref{eq3.8}), (\ref{eq3.9}) we consider the natural embedding
$$
\sum_1^{\infty} \gamma_k^0 \subset \partial N^4 (X_0^2) \subset N^4 (X_0^2) \subset \Delta_1^4 \, .
$$
Of course, the $\gamma_k^0$ bounds the $D^2 (\gamma_k^0)$ in $\Delta_1^4 - {\rm int} \, N^4 (X_0^2)$. But we can do much better than that. Consider
$$
\sum_1^{\infty} \gamma_k^0 \subset X_0^2 \subset \Delta_1^4 \, .
$$
We can (essentially) extend the $\underset{1}{\overset{\infty}{\sum}} \ \gamma_k^0$ to an embedded family of discs
\begin{equation}
\label{eq3.10}
\sum_1^{\infty} d_k^2 \longrightarrow \Delta_1^4 \, ,
\end{equation}
which touches the $X_0^2 \subset \Delta_1^4$ only along $\underset{1}{\overset{\infty}{\sum}} \ \gamma_k^0$ and which also {\ibf smears itself arbitrarily tightly close} to $X_0^2$ (contrary to the $\underset{1}{\overset{\infty}{\sum}} \ D^2 ( \gamma_k^0)$ which clearly does not). The ``essentially'' here stems from the fact that, in the case $\Delta^3 \times I$, the $d_k^2 = d^2 (\gamma_k^0)$ in (\ref{eq3.10}) are defined only for those normal $\gamma^0$'s not corresponding to the defects $\# \, \infty \, \# \, (S^2 \times D^2)$ of the infinite RED collar $X^4 - \Delta^4$. It will turn out, eventually, that there is no harm in this. The construction of (\ref{eq3.10}) makes an essential use of the RED 3-dimensional collapsing flow, about which not much can be said at the level of the present smallish paper. Finally, we still have to mention one of the important ingredients of the present approach: there is a certain {\ibf compatibility} property between the RED 2-dimensional and 3-dimensional collapsing flows; we will not make it explicit here, but just refer to it, when necessary.

\smallskip

Now, once we have (\ref{eq3.10}), we will forget about $X^2$ and only retain the $X_0^2$ from (\ref{eq3.8}). With the same (\ref{eq3.8}), let us notice the following feature of our present set-up. Let us consider any of the RED 1-handles of $\Delta^2$, namely the
$$
\sum_1^n R_i \subset \Gamma (1) \subset \Gamma (\infty) \cup \sum_1^{\infty} D^2 (C_j) \, , \quad \mbox{with} \ \Gamma (1) - \sum_1^n R_i = {\rm tree} \, .
$$
When one adds to any of these $R_i$ all the incoming trajectories of the RED 2-dimensional collapsing flow, then we get the object
\begin{equation}
\label{eq3.11}
\{ \mbox{extended cocore of} \ R_i \} \subset \Gamma (\infty) \cup \sum_1^{\infty} D^2 (C_j) \, ,
\end{equation}
which is an infinite PROPERLY embedded tree, which splits locally the target. Even better, we get this way a PROPERLY and properly embedded copy of $B^3 - \{$a tame Cantor set of $\partial B^3 \}$, which we denote just like in (\ref{eq3.11}), by
\begin{equation}
\label{eq3.12}
\{ \mbox{extended cocore} \ R_i \} \subset N^4 (\Gamma (\infty)) \cup \sum_1^{\infty} D^2 (C_j) \, .
\end{equation}
As a matter of terminology, by ``proper'' we mean boundary to boundary and interior to interior, while by ``PROPER'', in capital letters, we mean $f^{-1} ({\rm compact}) = {\rm compact}$.

\smallskip

Now, the same kind of construction as for (\ref{eq3.12}) also works perfectly well, in the following cases, for instance.

\medskip

i) Consider the set of the $b_i \in B \cap \Gamma (1)$. For obvious reasons, we have $\# \, B \cap \Gamma (1) \geq \# \, R \cap \Gamma (1) = n$ and we may as well assume that
\begin{equation}
\label{eq3.12.1}
P \underset{\rm def}{=} \, \# \, B \cap \Gamma (1) > \# \, R \cap \Gamma (1) = n \, .
\end{equation}

Any of these $b_i$'s also has an $\{$extended cocore $b_i\}$. Notice that the $b_i$'s are {\ibf not} exactly 1-handles of $\Delta^2$ (at least if we insist, as we normally do, to have a unique handle of index $0$). But the (\ref{eq3.12.1}) is a {\ibf disbalance} between Red and Blue for $\Delta^4 (\approx \Delta^2)$ which, later on, we will have to deal with.

\medskip

ii) Let $p \in X_0^2$ be any smooth point of some $D^2 (C)$ (but not of any $D^2 (\Gamma)$). Then $p$ also possesses a PROPERLY, but not quite properly, embedded $\{$exten\-ded cocore $(p)\}$, inside $X_0^2$ and/or $N^4 (X_0^2)$. Let us say that the embedding fails to be proper, along a small disc of $\partial B^3 - \{$the tame Cantor set$\}$. It is the RED 2-dimensional collapse $X_0^2 \to \Delta^2$ which creates, of course, these $\{$extended cocores$\}$, which are absent for $X^2$.

\medskip

Finally, one should notice two capital sins of our present set-up, as it is
\begin{eqnarray}
\label{eq3.13.1}
&&\mbox{Our $X_0^2$ (or $X^2$ itself, for that matter), possesses two not everywhere well-defined} \\
&&\mbox{2-dimensional collapsing flows, the RED and the BLUE ones.} \nonumber
\end{eqnarray}
But, generally speaking, and unfortunately for us as it turns out, the two kinds of trajectories cut through each other transversally. The global picture of the set $\{$RED 2-flow lines$\} \cup \{$BLUE 2-flow lines$\}$ is horribly complicated.
\begin{equation}
\label{eq3.13.2}
\mbox{There are no $\{$exterior cocore $q\}$ for points $q \in {\rm int} \, D^2 (\Gamma_i)$; and there is certainly no cure for this.}
\end{equation}
But then, on the road to those embedded exterior discs in cancelling position which we are eventually after (see the very beginning of this section), we most likely have some provisional substitute discs which we call now $\delta^2$, neither quite embedded nor quite exterior. This last thing means transversal contacts $\delta^2 \cap X_0^2 \subset \Delta_1^4$. These contacts may take the form
$$
q \in \delta^2 \cap D^2 (\Gamma_i) \subset \Delta_1^4 \, , \eqno (*)
$$
and here the lack of $\{$exterior cocore $q\}$ is, as we shall see, a serious potential danger. Hence, it would be very desirable to eliminate all the $q$'s like in ($*$) above. We will manage to do that, completely, in the case $\Delta^3 \times I$; see the stage III below. What we will manage to do in the Schoenflies case will be just to control the occurances ($*$) to a sufficient extent so that they become manageable. This is a good place to stress one basic difference between the two levels of our discussion. In the $\Delta^3 \times I$ case, the $\Delta^2$ is embedded in dimension three; in the Schoenflies case this is certainly not so.

\smallskip

Now, before we manage to start dealing with the two issues (\ref{eq3.13.1}), (\ref{eq3.13.2}), a lengthy prentice will have to be opened.

\bigskip

\noindent STAGE III, a prentice on compactifications. We consider now a compact bounded smooth 4-manifold, with only handles of index one and two. Call it $\Delta^4$; this could, of course, be the $N^4 (\Delta^2)$ from (\ref{eq3.9}), but we are supposed so be now at a higher level of generality. We consider, also, an open 4-manifold
\begin{equation}
\label{eq3.14}
X_0^4 = \Delta^4 + \{\mbox{handles of index one, called $h_i$, and handles of index two called} \ D^2 (C_j)\} \, . 
\end{equation}
Let us be slightly more specific about the way in which our handles are attached. We start by adding to $\Delta^4$ finitely many infinite trees, which we thicken in dimension four
\begin{equation}
\label{eq3.14.1}
\Delta^4 \cup \sum_1^f T_j \approx \Delta^4 \cup \sum_1^f N^4 (T_j) \, .
\end{equation}
Next, one adds the 1-handles to (\ref{eq3.14.1}) and, afterwards, the 2-handles too, to the resulting space.

\bigskip

\noindent {\bf Lemma 7.} -- {\it We assume now that the geometric intersection matrix $C \cdot h$ is of the easy id $+$ nilpotent type.

\smallskip

\noindent There exists then a natural smooth compactification $\hat X_0^4$ of $X_0^4$, with the following properties
\begin{equation}
\label{eq3.15.1}
\mbox{We have a diffeomorphism} \ \hat X_0^4 = \Delta^4 \cup \{ \mbox{collar} \ \partial \Delta^4 \times [0,1]\} \, .
\end{equation}
Here, in the RHS, the two pieces are glued along $\partial \Delta^4 \times \{ 0 \}$, making that $\partial \hat X_0^4 = \partial \Delta^4 \times \{ 1 \}$. [It should be stressed here that the formula above is just a diffeomorphism, and that the real life embedding $\Delta^4 \subset X_0^4 \subset \hat X_0^4$, where $\partial \Delta^4 \cap \partial \hat X_0^4 \ne \emptyset$, is not quite the one which it may suggest.]
\begin{equation}
\label{eq3.15.2}
\mbox{There is a compact subset $F \subset \partial \Delta^4 \times \{ 1 \}$ such that $X_0^4 = \hat X_0^4 - F$};
\end{equation}
the next point describes the structure of $F$.
\begin{eqnarray}
\label{eq3.15.3}
&&\mbox{There is a disjoint partition $F = F_0 \cup F_1$, with $\bar F_1 = F$ where $F_0$ is a tame Cantor set} \\
&&\mbox{and $F_1$ is the closed set of a $1$-dimensional {\ibf tame lamination} ${\mathcal L}$.} \nonumber
\end{eqnarray}
The $F_0$ is actually the sum of the end-point spaces of the $T_i$'s.}

\bigskip

The (easy) id $+$ nil form of the matrix $C \cdot h$ puts the $C_i$'s and $h_i$'s into a natural bijection, each bloc $h_i \cup D^2 (C_i)$ being a 4-ball. We have, with ``$C\ell$'' standing for closure
$$
C\ell \left( X_0^4 - \Delta^4 \cup \sum_1^f N^4 (T_j) \right) = \sum_i (h_i \cup D^2 (C_i))
$$
and, with this we define the following pair of smooth non-compact manifolds
\begin{equation}
\label{eq3.16}
(\mbox{LAVA}, \delta \, \mbox{LAVA}) 
= \left( \sum_i (h_i \cup D^2 (C_i)) , \partial \, \mbox{LAVA} \cap \partial \left( C\ell \left( X_0^4 - \Delta^4 \cup \sum_1^f N^4 (T_j) \right) \right) \right) \, .
\end{equation}
The LAVA and the $\delta \, \mbox{LAVA}$ are non compact bounded smooth manifolds, of dimensions four and three respectively. It is via the $\delta \, \mbox{LAVA} \subset \partial \, \mbox{LAVA}$, that our LAVA connects to the other world.

\smallskip

In the context of Lemma 7, the key fact is the following
\begin{eqnarray}
\label{eq3.16.1}
&&\mbox{The pair $(\mbox{LAVA}, \delta \, \mbox{LAVA})$ has the following {\ibf product property}. There is a diffeomorphism} \nonumber \\
&&(\mbox{LAVA} \cup \{\mbox{the lamination ${\mathcal L}$, added at infinity}\} , \delta \, \mbox{LAVA}) 
= (\delta \, \mbox{LAVA} \times [0,1] , \delta \, \mbox{LAVA} \times \{ 0 \}). 
\end{eqnarray}

It is important to notice here that our product property is a feature of the pair $(L , \delta L)$, and not just an absolute property of the space $L$ above.

\smallskip

The space of endpoints $e \left( \overset{f}{\underset{1}{\sum}} \, T_i \right) = \overset{f}{\underset{1}{\sum}} \, e (T_i)$ glues naturally both to $\biggl\{ \overset{f}{\underset{1}{\sum}} \, T_i$ and/or to $\overset{f}{\underset{1}{\sum}} \, N^4 (T_i)$ and $\Delta^4 \cup \overset{f}{\underset{1}{\sum}} \, N^4 (T_i) \biggl\}$, and then also to LAVA $\cup \, {\mathcal L}$. This allows us to compactify LAVA itself into the following object
$$
\mbox{LAVA}^{\wedge} = \mbox{LAVA} \cup {\mathcal L} \cup \sum_1^f e (T_i) \, .
$$

With all these things, the explicit definition of the $\hat X_0^4$ from Lemma 7 is, actually
\begin{eqnarray}
\label{eq3.16.2}
\hat X_0^4 = \Delta^4 \cup \sum_1^f N^4 (T_i) \cup \mbox{LAVA}^{\wedge}
\end{eqnarray}
where the second ``$\cup$'', i.e. the way in which $\mbox{LAVA}^{\wedge}$ glues to the rest, requires some specifications which we will not explain here. The product property which, appropriately stated is shared by $(\mbox{LAVA}^{\wedge} , \delta \, \mbox{LAVA})$ too, is {\ibf the} big virtue of lava, as far as we are concerned. With the explicit description of $\hat X_0^4$ given just above, the diffeomorphism (\ref{eq3.15.1}) is now a consequence of the product property. As already said before, (\ref{eq3.15.1}) is only a diffeomorphism, the real life formula behind it, from which also the correct embedding $\Delta^4 \subset \hat X_0^4$ is readable, is actually (\ref{eq3.16.2}).

\bigskip

\noindent {\bf Some comments.} A) When Lemma 7 is applied to something like the explicit context of stage II, then the $\{$extended cocore $R_i\} \subset N^4 \left( \Gamma (\infty) \cup \underset{1}{\overset{\infty}{\sum}} \, D^2 (C_j) \right)$ or the $\{$extended cocore $p\} \subset N^4 (X_0^2)$ get themselves compactified into objects which we will call
\begin{equation}
\label{eq3.16.3}
\{\mbox{extended cocore}\}^{\wedge} = \{\mbox{extended cocore}\} \cup \{\mbox{Cantor set}\} = B_{\rm smooth}^3 \, .
\end{equation}

\noindent B) Now, simple-mindedly, one might think that, locally, our $\mbox{LAVA} \, \cup \, \{{\mathcal L} \ \mbox{at infinity}\}$ is nothing but an object like
$$
\{\mbox{some extended cocore}\}^{\wedge} \times [0,1] \, .
$$
But this is certainly not so. It is indeed true that $(\mbox{LAVA}) \cup {\mathcal L}$ comes equipped with a surjection
$$
(\mbox{LAVA}) \cup {\mathcal L} \twoheadrightarrow \{\mbox{some highly non simply-connected train-track}\},
$$
but the fibers jump here quite wildly, as one moves around the train-track in question. There is not, even locally, a product.

\medskip

\noindent C) The little theory above can be generalized when we have handles of index one, two {\it and} three. Our $F = \underset{1}{\overset{f}{\sum}} \, e (T_i) \cup {\mathcal L}$ occurs then as accumulation points of a second, 2-dimensional lamination by planes, living at the infinity of a 4-dimensional object we call MAGMA. This may be useful in some situations, presumably.

\medskip

\noindent D) The compactification above is considerably more simple-minded than the so-called {\ibf strange compactification} from \cite{PoVI}. The only similarity is that a lamination occurs there too; but that one has both nontrivial holonomy and some nasty singularities. All such things are absent here.

\medskip

\noindent E) The following pair
\begin{equation}
\label{eq3.17}
\left( N^4 \left( \Gamma (\infty) \cup \sum_1^{\infty} D^2 (C_j) \right)^{\wedge} , \ \sum_1^n \{\mbox{extended cocore} \ R_i\}^{\wedge} \right) \, ,
\end{equation}
is a standard connected sum of $n$ copies of
$$
(S^1 \times B^3 , (*) \times B^3) \, .
$$
With this, $N^4 (X_0^2)^{\wedge}$ is now a smooth handlebody decomposition of $\Delta^4$, with $n$ red handles of index one and $\bar n$ handles of index two. Here, in the case $\Delta^3 \times I$ we have $\bar n = n$ and in the Schoenflies case $\bar n > n$. This kind of paradigm will be very useful, later on.

\bigskip

\noindent STAGE IV. We go back now to the two stings (\ref{eq3.13.1}), (\ref{eq3.13.2}), starting actually with (\ref{eq3.13.2}).

\smallskip

In the beginning, our way of proceeding will be purely 2-dimensional and abstract. By ``abstract'', we mean here that 4-dimensional incarnations or even maps into $4$-manifolds, are not yet considered. So, starting from
$$
\Gamma (1) \subset \Delta^2 \subset X_0^2 \subset X^2 \, ,
$$
let us define the following 2-complex which, at least for the time being will replace $X^2$ (think of it as being ``$X^2$ (old)'')
\begin{equation}
\label{eq3.18}
X^2 ({\rm new}) = X^2 \cup [(\Gamma (1) \times [0 \geq \xi_0 \geq -1]) \cup (\Delta^2 \times (\xi_0 = -1))] \, ,
\end{equation}
where the $\Gamma (1) \subset X^2$ is glued to $\Gamma (1) \times (\xi_0 = 0)$. We will decide, by decree, that from now on $\Delta^2 \times (\xi_0 = -1)$ is to be our
$$
\Delta^2 \approx \{\mbox{2-spine or 2-squeleton of $\Delta^4$}\} \, .
$$
The point here is that, with this
$$
\Delta^2 = \Delta^2 \times (\xi_0 = -1) \subset X^2 ({\rm new}) \, ,
$$
all our BLUE and RED features (at least the 2-dimensional ones, for the time being) are still with us. Red-wise, the $D^2 (\Gamma_i)$ are now the $D^2 (\Gamma_i) \times (\xi_0 = -1)$, the old $D^2 (\Gamma_i) \times (\xi_0 = 0)$ being declared $D^2 (\gamma^0)$'s. The idea is that, besides this, on $X^2 = X^2 ({\rm old}) \subset X^2 ({\rm new})$, the RED labels stay (essentially) put, and the RED collapse is proceeding according to the following general scheme
$$
X^2 ({\rm new}) \to (\Gamma (1) \times [0 \geq \xi_0 \geq -1]) \cup \Delta^2 \times (\xi_0 = -1) \to \Delta^2 \times (\xi_0 = -1) \, .
$$

When we move to BLUE, the general scheme is to start by the following decree
\begin{equation}
\label{eq3.18.1}
\mbox{The 2-cells} \ D^2 (\Gamma_{i}) \times (\xi_0 = -1) \ \mbox{are now $D^2 (\gamma^1)$'s.}
\end{equation}
With this one can crush all the $\Delta^2 \times [0 \geq \xi_0 \geq -1]$ part of $X^2 ({\rm new})$ onto $\Gamma (1) \times (\xi_0 = 0) = \Gamma (1) \subset X^2 ({\rm old})$ and then proceed on the remaining $X^2 ({\rm old})$, exactly as before, in the old case.

\smallskip

Notice that, with our new RED story as set up above, comes also an $X_0^2 ({\rm new}) \subset X^2 ({\rm new})$, defined on the same lines as (\ref{eq3.8}). The process $X^2 ({\rm old}) \Rightarrow X^2 ({\rm new})$ which we have just reviewed, is part of a bigger, still abstract transformation
\begin{equation}
\label{eq3.19}
X^2 ({\rm old}) \Longrightarrow X^2 ({\rm new}) = X^2 \Longrightarrow 2X^2 \supset 2X_0^2 \, ,
\end{equation}
to be defined, explicitly, later on. For expository purposes, we will give now, before we move to (\ref{eq3.19}), the gist of the way in which the transformation ${\rm old} \Rightarrow {\rm new}$, will be incarnated 4-dimensionally; only the easier case $\Delta^4 = \Delta^3 \times I$ will be discussed right now.

\smallskip

In the situation $\Delta^3 \times I$, at least, we have suggested in Figure 1 an immersion $f$ of $X^2 ({\rm new})$ into $X^4 = X^3 \times R$. This is the $X^4$ from (\ref{eq3.1}) and $f$ is part of the following diagram
$$
\xymatrix{
X^2 ({\rm new}) \ar[rr]^{f} \ar[drr]_{\pi \circ f}  &&X^4 = X^3 \times R \ar[rr]^{\qquad\pi_0} \ar[d]^{\pi} &&R \\ &&X^3
}
$$
which is actually a piece of the larger diagram (\ref{eq3.24}) below. Here are some explanations concerning Figure 1, from which the reader is supposed to read the diagram above. The regular oblique grid, which one should imagine infinite in all directions, is suggesting the embedding $X^2 ({\rm old}) \subset X^3 \times R$. We have tilted it so as to make the map $\pi \circ f$ generic. In the present $\Delta^3 \times I$ context, $\Delta^2 \subset X^2 ({\rm old})$ lives at level $t=0$.

\smallskip

\noindent It is suggested by the dotted fat line, along which $t=0$ and $\xi_0 = 0$ coincide. The plain fat line stands for $\pi \circ f (\Delta^2 \times (\xi_0 = -1))$, while the oblique thinly dotted lines stand for $\pi \circ f (\Gamma (1) \times [0 \geq \xi_0 \geq -1])$. Our $f$ is a generic immersion and some of the double points $fM^2 (f)$ are represented as fat points. The
$$
\xymatrix{
X^2 ({\rm new}) \ar[r]^{\ \ \pi \circ f} &X^3
}
$$
itself, is a singular $2$-dimensional polyhedron, with undrawable singularities, in the sense of \cite{Ga} and/or \cite{PoI}. Anyway, with all these things we may define $N^4 (X_0^2 ({\rm new}))$ as being the regular neighbourhood of the immersion $f \mid X_0^2 ({\rm new})$, happily disregarding the double points of $f$.  There is no harm with this, of course.

$$
\includegraphics[width=8cm]{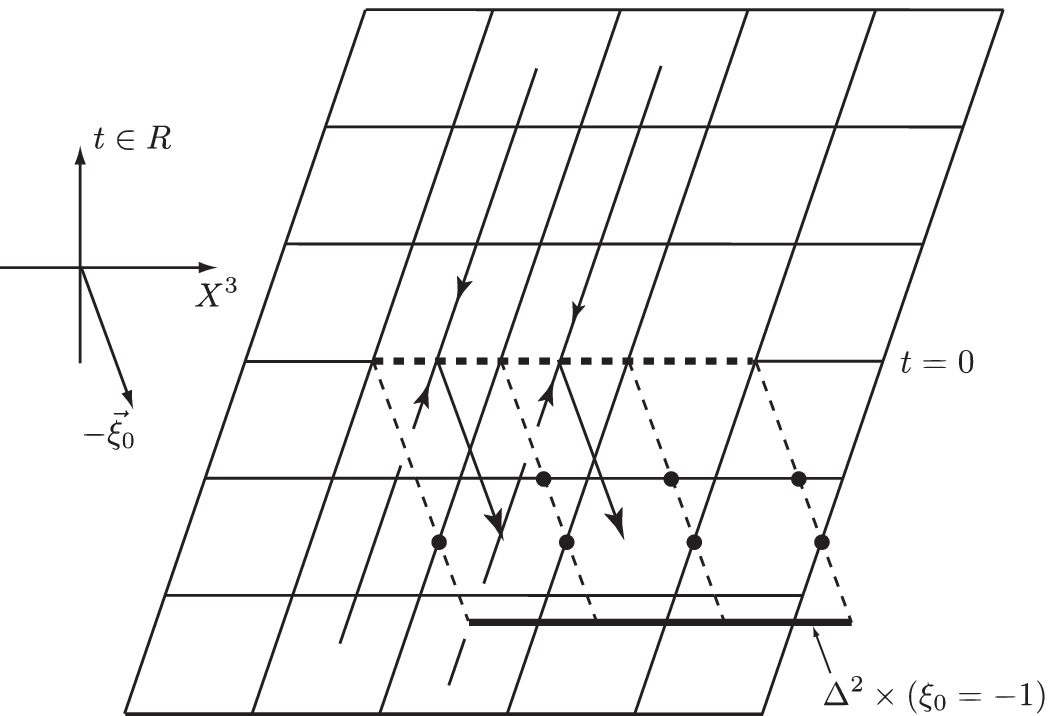}
$$
\begin{figure}[h]
\caption{$X^2 ({\rm new})$, in the context $\Delta^3 \times I$.}
\label{fig1}
\end{figure}

But Figure 1 is supposed to suggest other things too. Notice that in our drawing, the location of $f(\Gamma (1) \times [0 \geq \xi_0 \geq -1])$ breaks the symmetry between the past $(t < 0)$ and the future $(t > 0)$. The full incarnation of (\ref{eq3.19}) will break this symmetry even further. The point here is the following. Later on, curves like (\ref{eq3.4}), plus some others, will enter into a link projection (see (\ref{eq4.7}), in the next section) which will have to be changed into a link diagram. It so turns out that our geometrical set-up, when restricted just to $X_0^2 ({\rm old})$, makes that in the passage from the link projection to  the link diagram, there is a certain correlation between future and UP, and then also between past and DOWN. Of course, this convention, as such, could happily be reversed. But once it is there, then sending $[0 \geq \xi_0 \geq -1]$ to the past, like we did in Figure 1, will make that, in the context $\Delta^3 \times I$, the $\Delta^2 \times (\xi_0  = -1)$, actually the corresponding curves $\Gamma_i \times (\xi_0 = -1)$, will be constantly DOWN, in the link diagram. For reasons to be explained later, this will turn out to be very good for us. At this precise point, the Schoenflies case is different and also more difficult, as it will turn out.

\smallskip

At the level of Figure 1 we have suggested by arrows, very schematically, the RED 3-dimensional flow, crushing everything on $\Delta^2 \times (\xi_0 = -1)$. This flow is, in the context $\Delta^3 \times I$, perfectly ``smooth'', like in the old context, {\it except} for folding singularities at $\xi_0 = 0$. These do create serious technical problems and we will only be able to survive with them, at the price of a very heavy cure. In the Schoenflies case, the RED 3-dimensional collapsing flow comes with even more serious problems at $\xi_0 = -1$. We will come back to them later.

\smallskip

But right now, we go back to the specific issue (\ref{eq3.13.2}) for $\Delta^3 \times I$. It turns out that once the curves $\Gamma_{i \leq n} \times (\xi_0 = -1)$ are kept completely DOWN, from the viewpoint of our link diagram, then the $\Delta^2 \times (\xi_0 = -1)$ itself is kept disjoined from all the action to come, in particular from the kind of accidents which will be described by (\ref{eq3.33.2}) (and their yoked (\ref{eq3.33.1})), below. But since this is a bit too technical to be explained here, we will also give a more heuristical and intuitively easy, albeit less tight argument. In our $\Delta^3 \times I$ case, the $\Delta^2$ is embeddable in dimension three and, when we thicken in dimension four the immersion suggested by Figure 1 then, at the level of $N^4 (X_0^2 ({\rm new}))$ and/or $N^4 (X^2 ({\rm new}))$, we will get
$$
\Delta^2 \times (\xi_0 = -1) \subset \Delta^3 \subset \partial N^4 \, .
$$
This should make it plausible, at least, that it is now out of trouble. So, we go now to $\Delta^4$ Schoenflies, having in the back of our minds the same concern (\ref{eq3.13.2}). The first thing now, is to make sure that (\ref{eq3.18.1}) is {\ibf strictly} true. The issue here is the following. To begin with, once (\ref{eq3.18}) has the BLUE features, something like (\ref{eq3.18.1}) has to be there. But then, any further subdivision, and there will be many such, transforms, generally speaking, any $D^2 (\gamma^1)$ into a new, smaller $D^2 (\gamma^1)$, plus many $D^2 (\eta)$'s. Clearly this would spoil any initial strict (\ref{eq3.18.1}). In consequence, some hard work is necessary in order to maintain (\ref{eq3.18.1}) true, strictly. So, we assume this to be so, from now on, and let us see what it can do for us. Like in (\ref{eq3.5}) (and/or in (\ref{eq3.4})), in the new context (\ref{eq3.18}) we continue to have the following equality between infinite sets
\begin{equation}
\label{eq3.19.1}
\{ D^2 (\Gamma_i)\} + \{ D^2 (C_j)\} + \{ D^2 (\gamma_k^0)\} = \{ D^2 (\eta_{\ell}) \} + \{ D^2 (\gamma_m^1)\}
\end{equation}
and so, the (\ref{eq3.18.1}) implies that we also have
\begin{equation}
\label{eq3.19.2}
\{ D^2 (\eta_{\ell}) \} \subset \{ D^2 (C_j) \} + \{ D^2 (\gamma_k^0)\} \, ,
\end{equation}
while $\{ D^2 (\eta_{\ell}) \} \cap \{ D^2 (\Gamma_i)\} = \emptyset$. As we shall see later, this is crucial for handling the issue (\ref{eq3.13.2}) in the Schoenflies case. But, for the time being we leave it at that.

\smallskip

From now on, it will be understood that $X_0^2$, $X^2$ are the new ones (whenever the contrary is not explicitly said), and we move to the big issue (\ref{eq3.13.1}). In order to deal with it, we will introduce something like a much grander version of ((\ref{eq3.18}), the so-called {\ibf doubling process}. Beware that what follows next is, for the time being purely abstract, not yet incarnated 4-dimensionally.

\smallskip

Let $e \subset \Gamma(\infty)$ be an edge which contains an element $b_j \in B$; such an edge will be denoted, generically, by $e(b)$ (or sometimes, more specifically, $e(b_j)$). Any other edge, i.e. one containing either something in $R-B$ or nothing in $B \cup R$, will be denoted generically by $e(r)$. For the next purposes, we introduce three quantities
$$
0 < r \, (\mbox{for RED}) < \beta < b \, (\mbox{for BLUE}) \, ,
$$
with $\beta - r$ very small (compared to $b-r$).

\smallskip

For $b_i \in B$, we consider the curve $c_i (b) = \partial (e(b_i) \times [r,b])$, boundary of the 2-cell $D^2 (c_i(b)) = e(b_i) \times [r,b]$. Similarly, we define the 2-cell $D^2 (c(r)) = e(r) \times [r,b]$, cobounding $c(r) = \partial D^2 (c(r))$. With this, we consider now the following infinite 2-complex
\begin{equation}
\label{eq3.20}
2 X_0^2 = (X_0^2 \times r) \cup \{ \Gamma (\infty) \times [r,b] - \sum {\rm int} \, D^2 (c(b))\} \cup \left( \bigcup_{1}^{\infty} D^2 (\eta_{\ell}) \right) \times b \, ,
\end{equation}
a formula which begs for some explanations. The $X_0^2 \times r$ is glues to the middle term $\{ \ldots \}$ along $\Gamma (\infty) \times r$, while the middle term is then glued to
$$
X_b^2 \underset{\rm def}{=} \left(\bigcup_{1}^{\infty} D^2 (\eta_{\ell}) \right) \times b
$$
along $\Gamma (\infty) \times b$. Next, when in (\ref{eq3.20}) we delete ${\rm int} \, D^2 (c_i (b))$, it should be understood that a boundary collar, thicker than $\beta - r$ is left in place, so that the inclusion $c_i (b) \subset \partial (2 X_0^2)$ should make sense. Completely similarly, when the interiors of the 2-cells $D^2 (\gamma_k^0)$ (see (\ref{eq3.4}) and (\ref{eq3.5})) are deleted from $X^2 \approx X^2 \times r$, so as to get $X_0^2 \times r \subset 2 X_0^2$, again a collar is left in place so that, eventually, we should get
\begin{equation}
\label{eq3.20.1}
\partial(2 X_0^2) = \sum_k \gamma_k^0 + \sum_i c_i (b) \, .
\end{equation}
The 2-skeleton of $2 X_0^2$, which we denote by $2\Gamma (\infty)$ is
$$
2\Gamma (\infty) = (\Gamma (\infty) \times r) \cup (\Gamma_0 (\infty) \times [r,b]) \cup (\Gamma (\infty) \times b) \, ,
$$
where $\Gamma_0 (\infty) \subset \Gamma (\infty)$ is the 0-skeleton. Schematically speaking, $2 X_0^2$ consists of a red side $X_0^2 \times r$, a blue side $X_b^2$, plus some intermediary stuff which is essentially $\Gamma (\infty) \times [r,b]$, but with some deletions. On the same lines as in (\ref{eq3.20}), we introduce the following larger 2-complex
\begin{equation}
\label{eq3.21}
2 X^2 = (X_0^2 \times r) \cup (\Gamma (\infty) \times [r,b] ) \cup X_b^2 \, .
\end{equation}
In all this story, it should be understood that our object of interest, $\Delta^2 = \Delta^2 \times (\xi_0 = -1)$ lives now, naturally, in $X_0^2 \times r$; but the $D^2 (\eta)$'s have been transferred to the $b$-side.

\smallskip

In the next lemma, the RED 3-dimensional collapse is ignored.

\bigskip

\noindent {\bf Lemma 8.} -- {\it At the present abstract $2$-dimensional level of $2 X^2$, all the desirable RED and BLUE features are preserved. But moreover, for $2 X^2 \supset 2 X_0^2$ we also have the following}
\begin{equation}
\label{eq3.22}
\mbox{\it There are {\ibf no transversal intersections} between the RED flow-lines and the BLUE flow-lines.}
\end{equation}

\bigskip

The whole purpose of the doubling was exactly to get (\ref{eq3.22}). Here is also a sketch of proof for Lemma 8. To begin with on the $r$-side
$$
\Gamma (\infty) \times r \subset X_0^2 \times r
$$
we keep all the $R,B$ as they are, as well as the labels $D^2 (\Gamma)$, $D^2 (C)$. Of course, the $D^2 (\gamma^0)$'s are gone, but they have left a thin collar and the useful boundary piece $\gamma^0$, in their place. The edges $\Gamma_0 (\infty) \times [r,b]$ carry no $R \cup B$ labels, while each edge $e \times b \subset \Gamma (\infty) \times b$ will carry a (newly created) $R \cap B$. The set $\{ D^2 (\Gamma) \}$ does not change, but we will have extended sets of $C$, $\eta$, namely 
$$
\{\mbox{extended set of $C$'s}\} = \{ C \} + \{ c(r) \} + \{ \eta \times b \} \quad \mbox{and} \quad \{\mbox{extended set of $\eta$'s}\} = \{ \eta \times b \} + \{ c(b)\} + \{ c(r) \} \, .
$$
In the same vein, we have
$$
\{\mbox{extended set of $\gamma^0$'s}\} = \{ \gamma^0 \} + \{ c(b) \} \, , \ \{\mbox{extended set of $\gamma^1$'s}\} = \{\Gamma \} + \{ C \} \, .
$$
Notice that the first of these last two formulae, really makes $2 X_0^2$ be the analogue of $X_0^2$, after doubling, with $2 X^2$ in the role of $X^2$. On the $X_0^2 \times r$ side, the old RED geometric intersection matrix is kept as such. In the extended $C \cdot h$, which one can check to be of the easy id $+$ nil type, the $c(r)$ is dual to the corresponding $e(r) \times b$, while $\eta_i \times b$ is dual to $b_i \times b$. In the new BLUE geometric interaction matrix, which is again of the easy id $+$ nil type, the $\eta_i \times b$ dual to $b_i \times b$ and then, also, $c(b_i)$ is dual to $b_i \times r$ and $c(r)$ to $e(r) \times b$.

\smallskip

By $e(\ldots)$ we may mean the corresponding, newly created $R \cap B$.

\smallskip

The old matrix $\eta \cdot B$ finds itself transported now on the $b$-side, making (\ref{eq3.22}) possible.

\smallskip

The much larger $2X_0^2$ collapses now on our $\Delta^2$ and, with this, our discussion of Lemma 8 is finished.

\smallskip

So far, all this was purely abstract stuff. In order to incarnate it, geometrically, we start with the natural embedding $X_0^2 ({\rm old}) \subset X^3 \times R$, just like in the discussion coming with Figure 1. Next, one has to find a good way to extend it to a generic immersion
\begin{equation}
\label{eq3.24}
\xymatrix{
\Delta^2 \subset 2X_0^2 \ar[rr]^{f} \ar[drr]_{\pi \circ f}  &&X^4 = X^3 \times R \ar[rr]^{\qquad\pi_0} \ar[d]^{\pi} &&R \\ &&X^3
}
\end{equation}
We will have to come back to the (\ref{eq3.24}), but let us pretend it is with us now. What we do next, is the following item.
\begin{eqnarray}
\label{eq3.24.1}
&&\mbox{Completely disregarding the double points of the immersion $f$, one can get the regular} \\
&&\mbox{neighbourhood $N^4 (2X_0^2)$ which is supposed to contain the {\ibf correct} $\Delta^4 = N^4 (\Delta^2)$. } \nonumber
\end{eqnarray}

\noindent One can apply the compactification of Lemma 7 to $N^4 (2X_0^2)$ and get $N^4 (2X_0^2)^{\wedge}$ which comes then with a diffeomorphism
$$
N^4 (2X_0^2)^{\wedge} = N^4 (\Delta^2) \cup \{\mbox{collar}\} \, .
$$

Now, at the level of stage II we had the embedding (\ref{eq3.10}), which was smearing itself very tightly close to $X_0^2$. This last property means that we can carry now the $\underset{k}{\sum} \, d_k^2$ cobounding the $\underset{k}{\sum} \, \gamma_k^0 \subset \partial (2X_0^2)$ along, to our present context
\begin{equation}
\label{eq3.24.2}
\Delta^2 \subset X_0^2 \subset 2X_0^2 \hookrightarrow N^4 (2X_0^2)^{\wedge} \, .
\end{equation}
In other words, the $\sum \gamma_k^0$ extends now to an immersion which, for the better or for the worst replaces the by now deceased $\underset{k}{\sum} \, D^2 (\gamma_k^0)$
\begin{equation}
\label{eq3.24.3}
\sum_k d_k^2 \overset{\mathcal J}{\longrightarrow} N^4 (2X_0^2)^{\wedge} \, .
\end{equation}
For this generic immersion there are now both double points $M^2({\mathcal J}) \subset \underset{k}{\sum} \, d_k^2 \times \underset{k}{\sum} \, d_{\ell}^2 - \{\mbox{diagonal}\}$, at the source, inducing
\begin{equation}
\label{eq3.24.4}
x \in {\mathcal J} M^2 ({\mathcal J}) \subset N^4 (2X_0^2)^{\wedge} \, ,
\end{equation}
at the target and, also, transversal contacts
\begin{equation}
\label{eq3.24.5}
z \in {\rm Im} \,{\mathcal J} \cap 2X_0^2 \subset N^4 (2X_0^2)^{\wedge} \, .
\end{equation}
The contacts $z$ take the form
$$
z \in {\rm Im} \,{\mathcal J} \cap (X_0^2 \times r \cup \Gamma (\infty) \times [r,\beta]) \, ,
$$
and in the case $\Delta^3 \times I$ (but not necessarily so in the Schoenflies case) they are avoiding $\Delta^2$ altogether, like in the previous discussion around (\ref{eq3.18}). At this point, one should also keep in mind that (\ref{eq3.24.3}) is the living memory of the RED 3-dimensional collapsing flow which, at least in the context (\ref{eq3.18}) was still with us. It was already mentioned that in the Schoenflies context of (\ref{eq3.18}), this RED 3-dimensional flow did have complications at $\xi_0 = -1$. The offshot of these complications, {\it are} the transversal contacts, hinted at above,
\begin{equation}
\label{eq3.24.6}
z \in {\rm Im} \, {\mathcal J} \cap \Delta^2 (\mbox{Schoenflies}) \, .
\end{equation}
This is part of the issue (\ref{eq3.13.2}) for Schoenflies, the discussion of which is still not finished yet.

\smallskip

Instead of coming directly to grips with the full diagram (\ref{eq3.24}), a more indirect road will be profitable now. We will use the language of singular 2-dimensional polyhedra, their desingularizations and 4-dimensional thickenings, which is explained in great detail in \cite{Ga}; see also \cite{PoI} or \cite{PoII}. It is understood that, whatever $f$ itself may be, it is generic with respect to $\pi$, so that
\begin{equation}
\label{eq3.25}
\xymatrix{
2X_0^2  \ar[r]^{\ \ \pi \circ f} &X^3
}
\end{equation}
is a singular 2-dimensional polyhedron. The $\pi_0$-values in (\ref{eq3.24}) induce naturally a desingularization $\varphi$ for (\ref{eq3.25}). Essentially, this will be the following prescription
\begin{equation}
\label{eq3.26}
\mbox{$\varphi = S$ means high $\pi_0$ values and $\varphi = N$ means low $\pi_0$ values.}
\end{equation}
With this comes a 4-dimensional thickening $\Theta^4$ which is diffeomorphic to the regular neighbourhood of the immersion $f$, i.e.
\begin{equation}
\label{eq3.26.1}
N^4 (2X_0^2) = \Theta^4 (2X_0^2 , \varphi) \, ,
\end{equation}
an equality stemming directly from first principles. In the references \cite{Ga}, \cite{PoI}, \cite{PoII} it is explained, in detail, how to any pair (singular 2-dimensional polyhedra, desingularization), a canonical 4-dimensional thickening $\Theta^4 ( \ldots , \ldots)$ is attached. We will rather concentrate here on the following restriction of (\ref{eq3.25}), namely
\begin{equation}
\label{eq3.26.2}
\xymatrix{
(X_0^2 \times r) \cup \Gamma (\infty) \times [r,\beta] \ar[r]^{\qquad\qquad \pi \circ f} &X^3 \, ,
}
\end{equation}
where most of our head-aches will be concentrated. Here is a brief description of how (\ref{eq3.26.2}) goes. Start with the restriction of (\ref{eq3.26.2}) to $X_0^2 ({\rm old})$. There, our set-up is such that any of the undrawable singularities has to involve a purely spatial branch and a time-like branch. For these singularities, we will have, according to the case
\begin{equation}
\label{eq3.27}
\varphi (\mbox{future}) = S \, , \ \varphi (\mbox{past}) = N \, .
\end{equation}
Out of $X_0^2 ({\rm old})$ grow branches $\Gamma (\infty) \times [r,\beta]$ and $\Gamma (1) \times [0 \geq \xi_0 \geq -1]$, creating new singularities, but never do we find, simultaneously, at a given singularity, a branch $t < 0$ and a branch $\xi_0 < 0$. The correct set-up now, is the following extension of (\ref{eq3.27}), which also supersedes it, whenever that is the case
\begin{equation}
\label{eq3.28}
\varphi (\Gamma (1) \times [0 \geq \xi_0 \geq -1]) = N \, , \ \varphi (\Gamma (\infty) \times [r,\beta]) = S \, .
\end{equation}

We will not try to explain exactly here why this {\ibf is} the correct thing to do. Nor do we make explicit how to proceed at $\xi_0 = -1$. For the case $\Delta^3 \times I$ this is simple-minded, and rather clearly suggested by Figure 1. For the Schoenflies case, we also want the branches containing $+ \vec\xi_0$ to be always with $\varphi = S$, but this is now no longer hundred per cent automatic, and some work is needed at this particular point.

\medskip

\noindent {\bf Remark.} At any given singularity, exactly one of the four prescriptions in (\ref{eq3.27}) $+$ (\ref{eq3.28}) applies.

\smallskip

With all these things, the climax of the present stage IV, will be to replace the set-up (\ref{eq3.1}) by the following one
\begin{equation}
\label{eq3.29}
\Delta^4 = N^4 (\Delta^2) \subset N^4 (2X_0^2) \subset \Delta_1^4 \overset{\mathcal J}{\longleftarrow} \sum_k d_k^2 \, ,
\end{equation}
where we take as ambient space $\Delta_1^4$ a slightly larger copy of $N^4 (2X_0^2)^{\wedge}$, into which ${\mathcal J} d^2$ is pushed, rel its boundary, as much as it is possible.

\bigskip

\noindent STAGE V. ON THE WAY TO A SYSTEM OF EMBEDDED EXTERIOR DISCS, IN CANCELLING POSITION. The three adjectives just used will certainly not apply to the discs from Lemma 9 below, they remain for the time being at the level of a pipe dream.

\smallskip

We start now with the set $B$ {\ibf at level (\ref{eq3.18})} (and not at the fully doubled level (\ref{eq3.20}) or (\ref{eq3.21})), and we will organize it according to its natural BLUE order which comes from the easy id $+$ nil property of the matrix $\eta \cdot B$ (\ref{eq3.18}). Let us denote this by
$$
B = \{ b_1 , b_2 , b_3 , \ldots \} \, .
$$

Still at level (\ref{eq3.18}), with the same index $i$, we also have $e(b_i) \subset \eta_i \subset D^2 (\eta_i) \subset X_0^2$. Generically, the curve $\eta_i$ also has some other edges $e(b_j) \subset \eta_i$, with $j < i$. Now, in terms of the equality (\ref{eq3.19.1}), to the blue 2-cell $D^2 (\eta_i)$ corresponds a red cell which may be a $D^2 (\Gamma)$, a $D^2 (C)$ or a $D^2 (\gamma^0)$. Let us call this red cell $[D^2 (\eta_i)]$, but with the understanding that in the particular $D^2 (\gamma_k^0)$ case, then
$$
[D^2 (\eta_i)] = \{\mbox{the surviving collar of $D^2 (\gamma_k^0)$ in $X_0^2$}\} \, .
$$
In that case $\gamma_k^0$ occurs as an additional boundary piece of $[D^2 (\eta_i)]$. We move now to $2X_0^2$ and, at its level, for every $b_i \in B$ we define the following disc with holes which is properly embedded inside $2X_0^2$
\begin{eqnarray}
\label{eq3.30}
&&B_i^2 = \{ [D^2 (\eta_i)] \times r \} \cup \{ \eta_i \times [r,b] , \ \mbox{with} \ D^2 (c_i (b)) \ \mbox{and the} \ D^2 ((c_{j < i} (b)) \ \mbox{replaced by their} \\
&&\mbox{corresponding surviving collars, creating thus boundary pieces} \ c_i (b) , c_j (b) \ \mbox{for} \ B_i^2 \} \cup \{ D^2 (\eta_i) \times b\} . \nonumber
\end{eqnarray}
As a notational remark, do not mix up the $B_i^2$ from (\ref{eq3.30}) with the Blue set $B$, they are not at all the same thing.

\smallskip

Notice that
\begin{equation}
\label{eq3.31}
\partial B_i^2 = c_i (b) + \sum \{\mbox{some lower $c_{j<i} (b) \} + \{$possibly,} \ \gamma^0_k \} \, .
\end{equation}
Let us go back now, for a minute, to the immersion (\ref{eq3.24.3}), which also appears in (\ref{eq3.29}). The map ${\mathcal J}$ of the $d_k^2$ into $\Delta_1^4$ is guided by another map, which really does come from the RED 3-flow, namely
\begin{equation}
\label{eq3.32}
\sum_k d_k^2 \overset{F}{\longrightarrow} X_0^2 \subset 2X_0^2 \, .
\end{equation}
Think of $F$ as being more or less immersive too (but it certainly has folds) and with these things, the ${\mathcal J}$ (\ref{eq3.24.3}) is the lift of $F$ into $\Delta_1^4$.

\bigskip

\noindent {\bf Lemma 9.} -- 1) {\it For each curve $c_i (b)$ there is a disc $\delta_i^2$, with $\partial \delta_i^2 = c_i (b)$, coming with a map
\begin{equation}
\label{eq3.32.1}
\sum \delta_i^2 \overset{F}{\longrightarrow} 2X_0^2 \, .
\end{equation}
We construct} (\ref{eq3.32.1}) {\it by induction, along the natural BLUE order. So, for given $b_i$, assume that} (\ref{eq3.32.1}) $\vert \, (j < i)$ {\it is already well-defined. Then, for $c_i (b)$ we consider the $B_i^2$, with which $F\delta_i^2$ will start. Next, we have to fill in the missing boundary pieces occuring in} (\ref{eq3.31}). {\it We fill in every $c_{j<i} (b)$ with $F \delta_j^2$ and $\gamma_k^0$ (if it is there) with $Fd_k^2$.}

\smallskip

\noindent 2) {\it One can lift $F$ off $2X_0^2$ into an immersion which, essentially, extends {\rm (\ref{eq3.24.3})}, namely}
\begin{equation}
\label{eq3.32.2}
\sum \delta_i^2 \overset{\mathcal J}{\longrightarrow} \Delta_1^4 \, .
\end{equation}

\bigskip

Our ${\mathcal J}$, which rests on $N^4 (2X_0^2)^{\wedge} \subset {\rm int} \, \Delta_1^4$ (\ref{eq3.29}) exactly along $\sum c_i (b) \subset \partial N^4 (2X_0^2)^{\wedge}$ has, generally speaking, ACCIDENTS extending (\ref{eq3.24.4}), (\ref{eq3.24.5}), namely
\begin{equation}
\label{eq3.33.1}
\mbox{Double points} \ x \in {\mathcal J} M^2 ({\mathcal J}) \subset \Delta_1^4 \, , \ \mbox{and}
\end{equation}
\begin{equation}
\label{eq3.33.2}
\mbox{Transversal contacts} \ z \in {\mathcal J} \delta^2 \cap 2X_0^2 \subset \Delta_1^4 \, .
\end{equation}

\noindent A VERY IMPORTANT REMARK. In the case $\Delta^3 \times I$, we also have non-trivial $D^2 (\gamma_k^0)$ corresponding to the defects $\# \, \infty \, \# \, (S^2 \times D^2)$, and we will denote them, generically, by $D^2 (\gamma_{k(\beta)}^0)$. Obviously, our little story above would get into deep trouble every time we would find that
$$
[D^2 (\eta_i)] = D^2 (\gamma_{k(\beta)}^0) \, . \eqno (*_1)
$$

So, here is a hint how  we get rid of this problem. To begin with, not the full $\underset{1}{\overset{\infty}{\sum}} \, \delta_i$ will be actually needed, but only a very high finite truncation $\underset{1}{\overset{M}{\sum}} \, \delta_i^2$ of it. The quantity $M$ has to be large enough, if (\ref{eq3.32.2}) is to be good enough for our purposes. Next, at the level of (\ref{eq3.5}), there is a certain margin of flexibility for fixing the infinite subset $\sum D^2 (\gamma_{k(\beta)}^0) \subset \sum D^2 (\gamma_k^0)$. This turns also out to be an issue where there is no difference between $X^2 ({\rm old})$ and $X^2 ({\rm new})$. But the real point is now the following: we can fix $\sum D^2 (\gamma_{k(\beta)}^0)$, {\ibf after} having decided on the size of $M$. And then, we can also choose $\sum D^2 (\gamma_{k(\beta)}^0)$ sufficiently close to infinity so that, for $i \leq M$, the $(*_1)$ should not occur. There is no such problem in the Schoenflies case, of course. This ends our remark.

\medskip

Retain that each $\delta^2$ is made out of spare parts, possibly occuring with multiplicities, each of them being a $B_i^2$ or a $d_k^2$. As a result of the tri-dimensionality of $\Delta^3$ and of the passage old $\Rightarrow$ new, in the $\Delta^3 \times I$ context we will find that
\begin{equation}
\label{eq3.34}
{\mathcal J} \delta^2 \cap \Delta^2 = \emptyset \, .
\end{equation}
This is the complete happy end as far as the issue (\ref{eq3.13.2}) is concerned, in the case $\Delta^3 \times I$. As already noticed, on the other hand, in the Schoenflies case we do have
\begin{equation}
\label{eq3.34.S1}
{\mathcal J} d^2 \cap \Delta^2 \ne \emptyset \, .
\end{equation}
But then, once we have a {\ibf strict} (\ref{eq3.18.1}), and see here also the discussion coming with (\ref{eq3.19.1}), (\ref{eq3.19.2}), we will also find that
\begin{equation}
\label{eq3.34.S2}
{\mathcal J} B^2 \cap \Delta^2 = \emptyset \, .
\end{equation}
The two formulae above, in particular the last one, are what takes care of the issue (\ref{eq3.13.2}) in the Schoenflies case. Incidentally, it was quite an illumination for me when, during the Spring 2003 in Princeton, I realized that (\ref{eq3.18.1}), leading to (\ref{eq3.34.S2}), was the key to the until then locked door for Theorem 2.

\bigskip

\noindent {\bf Remark.} A priori, we might have tried to invoke (\ref{eq3.18.1}) for the case $\Delta^3 \times I$ too. For some technical reasons, connected to the 3-dimensional RED collapsing flow, we have chosen not to proceed that way.

\smallskip

As things stand right now, for our present $c_i (b) = \partial \delta_i^2$ when we move from the $B$ (\ref{eq3.18}) to the actual larger set of blue 1-handles $B$ of $2X_0^2$, then the contacts $c_i (b) \cdot B$ are exactly the following two
\begin{eqnarray}
\label{eq3.35}
&&c_i (b) \cdot (b_i \times r) = 1, \ \mbox{on the $X_0^2 \times r$ side, and} \\
&&c_i (b) \cdot (b_i \times b) = 1 , \mbox{on the $X_b^2$ side}. \nonumber
\end{eqnarray}
Let us also introduce the notation
\begin{equation}
\label{eq3.36}
\Delta^2 \cap B = \Gamma (1) \cap B = \{ b_{i_1} , b_{i_2} , \ldots , b_{i_P} \} \subset \{ b_1 , b_2 , \ldots , b_M \} \subset B \ \mbox{(\ref{eq3.18})} \, .
\end{equation}
The quantity $P$ here is the same one as in  (\ref{eq3.12.1}). In what follows next, in various successive steps, we will vastly change the system (\ref{eq3.32.2}). For these vastly transformed $\delta_i^2$'s, the boundary curve, which so far is the $c_i (b)$ above, will be denoted by $\eta_i ({\rm green}) = \partial \delta_i^2 \subset \partial N^4 (2X_0^2)^{\wedge}$. Also, with the large $B = B(2X_0^2)$, we will be very much focusing on the geometric intersection matrix
\begin{equation}
\label{eq3.37}
\eta ({\rm green}) \cdot B \, .
\end{equation}
Careful here, the $\delta^2$ is an {\ibf exterior} disk (or at least a candidate thereof). So, one should not mix up the {\ibf external} BLUE geometric intersection matrix (\ref{eq3.37}), with the related {\ibf internal} BLUE geometric intersection matrices $\eta \cdot B$ (\ref{eq3.18}) or $\eta \cdot B (2X_0^2)$.

\bigskip

\noindent {\bf Lemma 10.} -- {\it By sliding the external system of discs $\underset{1}{\overset{M}{\sum}} \, \delta_i^2$ over the internal $2$-handles contained in
$$
X_b^2 \cup \Gamma (\infty) \times [b,r] \subset 2X_0^2
$$
we can create a new system of external discs, which after further restriction from $M$ to $P$ {\rm (\ref{eq3.36})} we denote
\begin{equation}
\label{eq3.38}
\sum_{\ell = 1}^P \delta_{i_{\ell}}^2 \overset{\mathcal J}{\longrightarrow} \Delta_1^4 \, , \quad \partial \delta_{i_{\ell}}^2 = \eta_{i_{\ell}} ({\rm green}) \, ,
\end{equation}
which is such that the following {\ibf blue diagonality} condition should be satisfied:
\begin{equation}
\label{eq3.39}
\mbox{For $\alpha , \beta \leq P$ we have $\eta_{i_{\alpha}} ({\rm green}) \cdot b_{i_{\beta}} = \delta_{\alpha\beta}$, and} \ \sum_{\ell = 1}^P \eta_{i_{\ell}} ({\rm green}) \cdot (B(2X_0^2) - \Delta^2 \cap B (2X_0^2)) = 0 \, .
\end{equation}
}

\bigskip

This operation increases, a priori, the bag of ACCIDENTS to be considered, afterwards. Also, both very importantly and less trivially so than it might seem, there is no obstruction for performing this BLUE diagonalization. We will have to come back to this issue.

\smallskip

Once our diagonalization (\ref{eq3.39}) has been performed, we can afford to go to a simpler notation, namely
$$
\{ b_1 , b_2 , \ldots , b_P \} = \{ b_{i_1} , \ldots , b_{i_P} \} \ {\rm and} \ \eta_{i_j} ({\rm green}) = \eta_j ({\rm green}) \, , \ j \leq P \, .
$$

At this point, there are some very serious problems to be faced, which we list below.
\begin{eqnarray}
\label{eq3.40.1}
&&\mbox{We do have ACCIDENTS. Dealing with them is actually the hardest and longest part} \\
&&\mbox{of the proofs in PoV-B (\cite{PoV-B}).} \nonumber
\end{eqnarray}
But, for expositary purposes, we will pretend in the next Stage VI that the accidents have already been dealt with. In the next Section 4 some hints will be given concerning the operation of killing all the accidents, which in real life will have to preceed the $R/B$-balancing of Stage VI.
\begin{eqnarray}
\label{eq3.40.2}
&&\mbox{So, assume there are no accidents, and the (\ref{eq3.38}) is really a system of exterior discs.} \\
&&\mbox{We also want them to be in cancelling position with the 1-handles of $\Delta^4$. But then,} \nonumber \\
&&\mbox{these are the $R_1 + \cdots + R_n$ and {\ibf not} the $b_1 + \cdots + b_P$; the $\Gamma (1) - \underset{1}{\overset{P}{\sum}} \, b_i$ is {\ibf not} connected.} \nonumber
\end{eqnarray}
The reader may check that, in an ideal world where (\ref{eq3.40.1}) would have already been dealt with and where we would also have $P=n$, we would be done, by now. But $P > n$, in real life.

\smallskip

Now, (\ref{eq3.40.1}) has to be dealt with before we come to grips with (\ref{eq3.40.2}). We will show how to handle (\ref{eq3.40.2}) in the next Stage VI, but this will come then with another new, very serious problem, as we shall see. In a nutshell, this will be that
\begin{eqnarray}
\label{eq3.40.3}
&&\mbox{Once (\ref{eq3.40.2}) will have been dealt with, the blue diagonalization (\ref{eq3.39}) will no longer} \\
&&\mbox{be good enough, and another GRAND BLUE DIAGONALIZATION will be needed.} \nonumber
\end{eqnarray}

This will be one of the topics of the next Section 4.

\bigskip

\noindent STAGE VI. A CHANGE OF VIEWPOINT CONCERNING $\Delta^4$. During the change of viewpoint in question, the product structure $\Delta^3 \times I$ will get blurred too, but that is fine since by now, at this stage of the game, it has already served its purpose.

\smallskip

In our context used so far, we had, remember
$$
\Delta_1^4 = \{\mbox{ambient space} \ N^4 (2 X_0^2)^{\wedge} \cup ({\rm collar})\} = \Delta^4 \cup ({\rm collar}) \, ,
$$
where the first equality is a definition and the second one a diffeomorphism. Let us say that, up to now, our context has been
\begin{equation}
\label{eq3.41}
\Delta^4 \subset N^4 (2X_0^2) (\mbox{non-compact}) \subset {\rm int} \, \Delta_1^4 \subset \Delta_1^4 \, ,
\end{equation}
with $\Gamma (1) =\{$1-skeleton of $\Delta^4 \} \subset \Gamma (\infty) = \{$1-skeleton $\Gamma (\infty) \times r$ of $X_0^2 \times r \} \subset 2\, \Gamma (\infty) = \{$1-skeleton of $2X_0^2 \}$. 

\smallskip

With this, as things stand now, we also have $\Gamma (1) \cap B = \underset{1}{\overset{P}{\sum}} \, b_i$, $\Gamma (1) \cap R = \underset{1}{\overset{n}{\sum}} \, R_i$, where $P \geq n$ and where also, since the case $P=n$ is easier, it will be assumed that $P > n$.

\medskip

\noindent {\bf The $R/B$ balancing Lemma 11.} -- {\it Staying all the time embedded inside the ambient space $\Delta_1^4$ and also keeping $\Delta^4$ fixed, we can submit the $N^4 (2X_0^2)$ in {\rm (\ref{eq3.41})} to the following kind of compact changes, localized inside $X_0^2 \times r$
\begin{equation}
\label{eq3.42.1}
\mbox{Pick up a certain well-chosen family of $1$-handles $\{ y_1 , y_2 , \ldots , y_{P-n}\} \subset \Gamma (\infty) \cap h - B \subset R-B$} \, ,
\end{equation}
the dual $C$-curves of which we will denote by $C(1), C(2) , \ldots , C(P-n)$.
\begin{eqnarray}
\label{eq3.42.2}
&&\mbox{We will perform now embedded $1$-handle slidings, dragging along the corresponding} \\
&&\mbox{$2$-handles, at the level of $N^4 (\Gamma (\infty)) \subset N^4 (2\Gamma (\infty))$. We will slide, in succession, each of the} \nonumber \\
&&\mbox{$y_1 , y_2 , \ldots , y_{P-n}$ over a second family of well-chosen elements $x \in h-B \subset R-B$ which are} \nonumber \\
&&\mbox{always such that, in the natural RED order of $C \cdot h$ {\rm (\ref{eq3.18})}, and hence of $C \cdot h \, (2X_0^2)$ too, we} \nonumber \\
&&\mbox{should have the following inequality, every time $y$ slides over $x$} \nonumber
\end{eqnarray}
$$
x < \{\mbox{the $y$ which slides}\}
$$
\begin{eqnarray}
\label{eq3.42.3}
&&\mbox{At the end of the sliding move, we have a $\Gamma (\infty) \subset 2 \Gamma (\infty)$ which has changed (but we do} \\
&&\mbox{not bother to denote these objects differently). The point is that $\Gamma (1)$ has been replaced} \nonumber \\
&&\mbox{by a $\Gamma (3) \subset \{\mbox{new} \ \Gamma (\infty)\}$, which is {\ibf well-balanced}, in the sense that the $\Gamma (3) \cap B = \underset{1}{\overset{P}{\sum}} \, b_i$} \nonumber \\
&&\mbox{(as before) and the $\Gamma (3) \cap R = \underset{1}{\overset{n}{\sum}} \, R_i + \underset{1}{\overset{P-n}{\sum}} \, y_j$ are now two sets of the same cardinality,} \nonumber \\
&&\mbox{with $\Gamma (3) - B$, $\Gamma (3) - R$ being, both, trees.} \nonumber
\end{eqnarray}
\begin{equation}
\label{eq3.42.4}
\mbox{We also have, more globally, that $2\Gamma (\infty) - R$ and $2\Gamma(\infty)-B$ are trees.}
\end{equation}
}

Up to a diffeomorphism which does not budge $\Delta^4 , \Delta_1^4$, the sequence (\ref{eq3.41}) does not feel the effect of the balancing process.

\smallskip

Ideally, we would be quite happy if we could argue, at this point, as follows. Decide that the 1-skeleton of $\Delta^4$ is now $\Gamma (3)$, which comes equipped with two sets of 1-handles among which we might chose one. These are the BLUE $b_1 , b_2 , \ldots , b_P$ and the RED $R_1 , R_2 , \ldots , R_n$, $R_{n+1} = y_1 , \ldots$, $R_P = y_{P-n}$. The $\Delta^4$ should have now a handlebody decomposition which we will call ``ideal'', with 2-handles
\begin{equation}
\label{eq3.43}
\{ D^2 (\Gamma_i) \} \, , \ D^2 (C(1)) , \ldots , D^2 (C (P-n)) \, .
\end{equation}
Notice that not only have the $y_1 , \ldots , y_{P-n}$ been {\ibf promoted} as 1-handles of $\Delta^4$, but also their dual $D^2 (C(1)) , \ldots$, $D^2 (C(P-n))$ as 2-handles of $\Delta^4$ too. Since, clearly
$$
\{ C(1) , \ldots , C(P-n) \} \cdot \{ y_1 , \ldots , y_{P-n}\} = {\rm id} + {\rm nilpotent} \, ,
$$
formally at least we are O.K. Also, provided that the ACCIDENTS (\ref{eq3.33.1}), (\ref{eq3.33.2}) have been killed, Lemma 10 would then provide us with exterior, embedded 2-handles $\delta^2$, in cancelling position with the blue 1-handles. This little ideal scenario has a very serious flaw, which is the following
\begin{equation}
\label{eq3.44}
\mbox{The 2-handles $D^2 (C(1)) , \ldots , D^2 (C(P-n))$ are, generally speaking, {\ibf not} directly attached to $\Gamma (3)$.}
\end{equation}
Moreover, the (\ref{eq3.44}) is {\ibf irreparable} because of the following item. Notice, to begin with, that it is not hard to concoct a RED diagonalization process which would let some internal RED 2-handles $D^2 (C_a)$ slide over lower $D^2 (C_b)$, with $b < a$ in the natural red order, so that we should achieve
$$
C(i) \subset \partial N^4 (\Gamma (3)) \, , \eqno (*)
$$
But here comes a fact, which will be explained in Section 4
\begin{eqnarray}
\label{eq3.44.1}
&&\mbox{Contrary to the BLUE diagonalization which has led to (\ref{eq3.39}), the kind of RED} \\
&&\mbox{diagonalization leading to $(*)$, which was envisioned above, is {\ibf forbidden}.} \nonumber
\end{eqnarray}
But before we can explain what we will do now, in order to cope with these issues, we have to be more precise about our notations concerning the cardinalities of the handles of $\Delta^2$. The two $n = {\rm card} (\Gamma (1) \cap R)$, $P = {\rm card} (\Gamma (1) \cap B) = {\rm card} (\Gamma (3) \cap B)$ should be unambiguously clear and we use for them the same notations in both the contexts $\Delta^3 \times I$ and $\Delta^4$ Schoenflies. The cardinality of the set $\{ D^2 (\Gamma_i)\}$ used so far, is the same $n$ as above, in the case $\Delta^3 \times I$ but then also, it is some $\bar n > n$ in the Schoenflies case.

\smallskip

So, once the ideal scenario has collapsed, here is what we will do, in the real world. Consider, for the time being the {\ibf purely abstract promotion}, for $1 \leq i \leq P-n$
\begin{equation}
\label{eq3.45}
y_i \Longrightarrow R_{n+i} \ (\mbox{1-handle of} \ \Delta^4) \, , \ {\rm and}
\end{equation}
$C(i) \Longrightarrow \Gamma_{n+i}$ (case $\Delta^3 \times I)$, respectively $C(i) \Longrightarrow \Gamma_{\bar n + i}$ (case $\Delta^4$ Schoenflies); in both of these last two formulae, $\Gamma$ is the same physical curve as the $C$, but considered now as an internal attaching curve of $\Delta^4$.

\smallskip

With this promotion, which so far is barely more than a notational device, $\Delta^4$ is endowed now, abstractly speaking with $P$ handles of index one (a RED collection and also a BLUE one) and with the 2-handles
$$
\sum_1^P D^2 (\Gamma_i) \ \mbox{in the case} \ \Delta^3 \times I , \ \mbox{respectively} \ \sum_{1}^{P+(\bar n - n)} D^2 (\Gamma_i) \ \mbox{in the Schoenflies case.}
$$
Keep in mind that this is only abstract, so far, in the sense that for our physical $\Delta^4 = N^4 (\Delta^2)$, as such, no bonafide handlebody decomposition on the lines above is available. On the other hand, in this abstract context, when we look at the geometric intersection matrices all the desirable RED and BLUE features continue to be satisfied, provided that we also accompany the promotion (\ref{eq3.45}) by the following other related transformation, concerning now the whole of $N^4 (2\Gamma (\infty))$, after the $R/B$-balancing and promotion.
\begin{equation}
\label{eq3.46}
\mbox{By decree, the new families $h(2X_0^2)$ and $C(2X_0^2)$ are now the very slightly reduced}
\end{equation}
$$
h - \sum_{1}^{P-n} y_i \, , \ C- \sum_1^{P-n} C(i) \, , \ \mbox{respectively}.
$$
By the same decree we exclude from LAVA the $P-n$ copies of $B^4$ consisting of the $h_1 \cup D^2 (C(1)) , \ldots , h_{P-n} \cup D^2 (C(P-n))$. The new, slightly reduced lava, call it again $({\rm LAVA} , \delta \, {\rm LAVA})$, continues to have the product property. Moreover, via its $\delta \, {\rm LAVA}$, this new LAVA glues to
\begin{equation}
\label{eq3.47}
N^4 (2\Gamma (\infty) \ (\mbox{after $R/B$ balancing}) - h \ \mbox{(\ref{eq3.46}) (after promotion})) \supset N^4 (\Gamma (3)) \, .
\end{equation}
The next lemma, which should be compared to the comment E) at the end of Stage III, is our way to meet the difficulty (\ref{eq3.44}).

\bigskip

\noindent {\bf Lemma 12.} -- 1) {\it In the context of {\rm (\ref{eq3.47})} above, we introduce the following large smooth compact $4$-manifold
\begin{eqnarray}
\label{eq3.48}
\bar N^4 (\Gamma (3)) &= &\{ N^4 (2\Gamma (\infty) (\mbox{after $R/B$ balancing})) - h \mbox{\rm (\ref{eq3.46})}\} \cup ({\rm LAVA}^{\wedge}) \\
&= &(N^4 (2\Gamma (\infty)(\mbox{after balancing}) \cup \sum_{1}^{\infty} D^2 (C_i) (\mbox{after promotion}))^{\wedge} \, , \nonumber
\end{eqnarray}
where in the second term, the two pieces are glued along $\delta \, {\rm LAVA}$.

\smallskip

This $\bar N^4 (\Gamma (3))$ is a smooth compact handlebody of genus $P$, by which we mean, it is a $P \, \# \, (S^1 \times B^3)$.}

\medskip

2) {\it The $\bar N^4 (\Gamma (3))$ also comes equipped with a properly embedded system of $P$ $3$-balls ($=$ $1$-handle cocores), namely the
\begin{equation}
\label{eq3.49}
\sum_1^P \{\mbox{extended core} \, b_i \}^{\wedge} \subset \bar N^4 (\Gamma (3)) \, ,
\end{equation}
and with this, the pair defined by {\rm (\ref{eq3.49})} is {\ibf standard}.}

\medskip

3) {\it With our abstract promotion presented above in mind, we introduce now the following quantity
$$
\bar P \underset{\rm def}{=} P \, (\mbox{in the case} \ \Delta^3 \times I) , \bar P \underset{\rm def}{=} P + (\bar n - n) (\mbox{in the $\Delta^4$ Schoenflies case}).
$$
With this, the $2$-handles $\underset{1}{\overset{\bar P}{\sum}} \ D^2 (\Gamma_i)$ (which include now the promoted $D^2 (C(1)) , \ldots , D^2 (C(P-n))$, are quite naturally, directly attached to $\bar N^4 (\Gamma (3))$. Also, there is a diffeomorphism
\begin{equation}
\label{eq3.50}
\Delta^4 \underset{\rm DIFF}{=} \bar N^4 (\Gamma (3)) + \sum_1^{\bar P} D^2 (\Gamma_i) \, ,
\end{equation}
where, remember, $\Delta^4$ is here our $N^4 (\Delta^2)$ and, in the rest of the paper, its incarnation will be the RHS of the formula} (\ref{eq3.50}).

\bigskip

In other words, $\Delta^4$ has now a smooth handlebody decomposition with $P$ 1-handles $\underset{1}{\overset{P}{\sum}} \{\mbox{extended cocore} \, b_i\}^{\wedge}$ and with 2-handles $\underset{1}{\overset{\bar P}{\sum}} \, D^2 (\Gamma_j)$ (taking the promotion here into account).

\smallskip

Here are some comments. To begin with, the sets $(h,C)$ occuring here are the ones from after doubling, slightly diminished by the promotion. Also, the RED 1-handles $\underset{1}{\overset{P}{\sum}} R_i$ of $\Delta^4$ can and will be forgotten.

\smallskip

Assuming the accidents already killed, our latest transformation of (\ref{eq3.1}), after the (\ref{eq3.9}) and (\ref{eq3.29}), is now the following, with the same $\Delta_1^4$ as in (\ref{eq3.29})
\begin{equation}
\label{eq3.51}
\bar N^4 (\Gamma (3)) \cup \sum_1^{\bar P} D^2 (\Gamma_i) \subset \Delta_1^4 \overset{\mathcal J}{\longleftarrow} \sum_1^P \delta_i^2 \, ,
\end{equation}
where ${\mathcal J}$ is an embedding into
$$
\Delta_1^4 - {\rm int} (\bar N^4 (\Gamma (3)) \cup \sum_1^{\bar P} D^2 (\Gamma_i)) \, ,
$$
and where, isotopically speaking
$$
\bar N^4 (\Gamma (3)) + \sum_1^{\bar P} D^2 (\Gamma_i) = N^2 (2X_0^2)^{\wedge} = \{\mbox{closure of} \ N (2X_0^2) \subset \Delta_1^4 \} \, .
$$

Moreover, as a consequence of the BLUE diagonalization (\ref{eq3.39}), for the $\eta_i ({\rm green}) = \partial \delta_i^2 \subset \partial (\bar N^4 (\Gamma (3)) \cup \underset{1}{\overset{\bar P}{\sum}} \, D^2 (\Gamma_i))$ we have now
\begin{equation}
\label{eq3.52}
\eta_i ({\rm green}) \cdot b_j = \delta_{ij} \ {\rm if} \ i,j \leq P \ {\rm and} \ \sum_1^P \eta_i ({\rm green}) \cdot \left( B(2X_0^2) - \sum_1^P b_i \right) = 0 \, .
\end{equation}
But, the new problem which has been created now, is that the 1-handles of $\bar N^4 (\Gamma (3)) \cup \underset{1}{\overset{\bar P}{\sum}} \, D^2 (\Gamma_i)$ are not exactly the BLUE $\underset{1}{\overset{P}{\sum}} \, b_i$, but the more exotic $\underset{1}{\overset{P}{\sum}} \, \{\mbox{extended cocore} \ b_i\}^{\wedge}$. This is the difficulty mentioned in (\ref{eq3.40.3}). What we find now is exactly the following
\begin{eqnarray}
\label{eq3.53}
&&\eta_i ({\rm green}) \cdot \{\mbox{extended cocore} \, b_j\}^{\wedge} = \delta_{ij} + \{\mbox{an additional, call it off-diagonal term coming} \\
&&\mbox{from those contacts} \ \eta_i ({\rm green}) \cdot h_k \ {\rm with} \ h_k \in h-B \ {\rm and} \ h_k \subset \{\mbox{extended cocore} \ b_j \}\} \, . \nonumber
\end{eqnarray}
There are finitely many $h_k$'s involved in (\ref{eq3.53}), all living inside $X_0^2 \times r \subset 2X_0^2$.

\smallskip

At the point which we have reached now, there are still two main obstacles between us and what we want to achieve, namely

\smallskip

1) We still have to get rid of the accidents.

\smallskip

2) After that has been done, we still have to achieve the GRAND BLUE DIAGONALIZATION, by which we mean the following
$$
\eta_i ({\rm green}) \cdot \{\mbox{extended cocore} \ b_j \}\}^{\wedge} = \delta_{ij} \, ,
$$
which of course, is equivalent to
$$
\eta_i ({\rm green}) \cdot \{\mbox{extended cocore} \ b_j \}\} = \delta_{ij} \, .
$$

The next section is entirely devoted to these two pending issues. But since, in real life, this story is considerably more technical than what has been going on so far, the exposition will be even more sketchy and impressionistic.

\section{Some additional technicalities}\label{sec4}
\setcounter{equation}{0}

We consider now the stage when the little blue diagonalization (\ref{eq3.52}) has been already achieved, but all the accidents (\ref{eq3.33.1}), (\ref{eq3.33.2}) of
\begin{equation}
\label{eq4.1}
\sum_1^P \delta_i^2 \overset{\mathcal J}{\longrightarrow} \Delta_1^4 \, , \quad \partial \delta_i^2 = \eta_i ({\rm green})
\end{equation}
are still with us. Normally, the double points $x$ (\ref{eq3.33.1}) and the transversal contacts $z$ (\ref{eq3.33.2}) come {\ibf yoked} together, and here is a toy-model for a typical system of yoked accidents. In some coordinate neighbourhood $U \subset \Delta_1^4$, the set $2X_0^2 \cap U$ consists of two transversal planes $Q_1 , Q_2$ with $Q_1 \cap Q_2 = P$, while ${\rm Im} \, {\mathcal J} \cap U$ consists of two smooth branches $A_1 , A_2$, parallel copies of $Q_1 , Q_2$ respectively, coming with
\begin{equation}
\label{eq4.2}
A_1 \cap Q_2 = z_2 \, , \ A_2 \cap Q_1 = z_1 \, , \ A_1 \cap A_2 = x \, ;
\end{equation}
the notation of (\ref{eq3.33.1}), (\ref{eq3.33.2}) are being used here. Let us start with the following remark.
\begin{equation}
\label{eq4.3}
\mbox{Assume that $z_1$ possesses an $\{$extended cocore $z_1 \} \subset 2X_0^2$.}
\end{equation}
Then, we can push $A_2$ over the compactified $\{$extended cocore $z_1 \}^{\wedge}$, like it is suggested to do, with dotted lines, in Figure~2, and get rid of $z_1$. This process does not change $\eta ({\rm green})$ at all, and we will certainly make use of it. Also, according to conditions, $x$ might be destroyed together with $z_1$ too.

\smallskip

But, at this point, it is not hard to see, and the reader should certainly try to figure this out alone, that even if both $\{$extended cocore $z_1 \}$ and $\{$extended cocore $z_2 \}$ are present (which will be the case, most of the times), we still {\ibf cannot} use the mechanism (\ref{eq4.3}) in order to completely destroy the yoked system $(z_1 , z_2 , x)$.

\smallskip

What one should retain is that even with the extended cocore mechanism available at $z_1$, something else is still necessary for dealing with $z_2$. This finishes the discussion of the toy-model, and we go back now to the singular 2-dimensional polyhedron (\ref{eq3.26.2}) with its desingularization (\ref{eq3.26}) (see here (\ref{eq3.27}) and (\ref{eq3.28}) too). Notice that the quantity $\beta$ in $(X_0^2 \times r) \cup (\Gamma (\infty) \times [r,\beta])$ is close enough to $r$ so as not to see the deletions $2X^2 - 2X_0^2$.

\smallskip

The 1-skeleton of $(X_0^2 \times r) \cup (\Gamma (\infty) \times [r,\beta])$ is the following object
\begin{equation}
\label{eq4.4}
\Gamma (\infty) = \Gamma (\infty) \times r \, , \ \mbox{with a little are} \ P \times [r,\beta] \ \mbox{sticking out of each vertex} \ P \in \Gamma (\infty) \, .
\end{equation}

With this, we will review now the construction of
\begin{equation}
\label{eq4.5}
\Theta^4 ((X_0^2 \times r) \cup (\Gamma (\infty) \times [r,\beta]) , \varphi) \subset \Theta^4 (2X_0^2 , \varphi) = N^4 (2X_0^2) \, ,
\end{equation}
which, as far as accidents go, is the most important part, since it houses $Fd^2$. We will refer now to the procedures explained in \cite{Ga} (see also \cite{PoI}, \cite{PoII}) and use the singular 3-dimensional version, rather than the 2-dimensional one. We may assume that the restriction of the map $\pi \circ f$ (\ref{eq3.26.2}) to the 1-skeleton (\ref{eq4.4}) is an embedding. Its regular neighbourhood is an infinite solid torus $N^3 (\Gamma (\infty))$, coming with some additional structures. We will introduce the notation
\begin{equation}
\label{eq4.6}
\Sigma_{\infty}^2 = \partial N^3 (\Gamma (\infty)) \, ,
\end{equation}
and this infinite open surface $\Sigma_{\infty}^2$ comes with a PROPERLY embedded system of small disks, which we call generically $\beta$. This is the trace of the $\underset{P}{\sum} \, P \times [r,\beta]$. Next, $\Sigma_{\infty}^2$ comes equipped with an infinite {\ibf link projection}
\begin{equation}
\label{eq4.7}
\sum_{1}^{\bar P} \Gamma_i + \sum_1^{\infty} C_j + \sum_1^{\infty} \gamma_k^0 + \sum_1^{\infty} [c_{\ell} (b \, {\rm or} \, r)] \overset{j}{\longrightarrow} \Sigma_{\infty}^2 \, ,
\end{equation}
with the following specifications. Each $[c_{\ell} (b \, {\rm or} \, r)]$ is a piece of the corresponding $c_{\ell} (b \, {\rm or} \, r)$, essentially
$$
\mbox{``} c_{\ell} (b \, {\rm or} \, r) \cap [(X_0^2 \times r) \cup (\Gamma (\infty) \times [r,\beta])]\mbox{''.}
$$
Concretely, $j [c_{\ell} (b \, {\rm or} \, r)]$ is an arc hooked at two spots $\beta$. The $j$ is a generic immersion (in particular it has no triple points), and it injects on each connected component of the L.H.S. of (\ref{eq4.7}). Never mind here that the interiors of $D^2 (\gamma_k^0)$, $D^2 (c_{\ell} (b))$ have been deleted, their very useful surviving collars and boundary pieces are still with us. We consider the double points $s \in jM^2 (j) \subset \Sigma_{\infty}^2$, and the main facts are here the following
\begin{eqnarray}
\label{eq4.8}
&&\mbox{There is a canonical bijection} \\ 
&&\mbox{$\{$the undrawable singularities of the singular $2$-dimensional polyhedron (\ref{eq3.26.2})$\} \approx jM^2 (j)$} \, .  \nonumber
\end{eqnarray}
\begin{equation}
\label{eq4.9}
\mbox{Each $s \in jM^2(j)$ lives, inside $\Sigma_{\infty}^2 = \partial N^3 (\Gamma (\infty))$, close to some {\it canonically} attached vertex $P \in \Gamma (\infty)$.}
\end{equation}
We can consider, inside $X^3$, source of $f$ (\ref{eq3.26.2}) a disjoined system of 3-balls $B^3 (P)$ each centered at $f(P)$, with radii much thicker than the width of $N^3 (\Gamma (\infty))$ and, with this, all the interesting part of the link projection (\ref{eq4.7}) lives inside
$$
\sum_P B^3 (P) \cap \Sigma_{\infty}^2 \, .
$$
Now, what we know from first principles, is that the desingularization $\varphi$ (\ref{eq3.26}) of (\ref{eq3.26.2}), gives a recipee for undoing the double points $s \in jM^2 (j)$. At each $s$, the ${\rm Im} \, j$ has two branches and, keeping in mind (\ref{eq4.8}), we pull UP, towards the observer, the branch coming with $\varphi = S$ and then, accordingly to this, we push DOWN, the one with $\varphi = N$. Our set-up in (\ref{eq3.28}) makes that
\begin{equation}
\label{eq4.9bis}
\varphi [c_{\ell} (b \, {\rm or} \, r)] = S
\end{equation}
making that, whenever this makes sense, ``UP'' looks towards $b$ and ``DOWN'' towards $r$, with $b,r $ standing for blue and red, respectively, too. Keep in mind that all these are mere useful conventions.

\smallskip

We go now 4-dimensional and for this, we start by changing $N^3 (\Gamma (\infty))$ into $N^4 (\Gamma (\infty)) = N^3 (\Gamma (\infty)) \times [0,1]$, with $\partial N^4 (\Gamma (\infty))$ equal to the double of $N^3 (\Gamma (\infty))$. Very explicitely, we have now a {\ibf splitting}
\begin{equation}
\label{eq4.10}
\partial N^4 (\Gamma (\infty)) = \partial^- N^4 (\Gamma (\infty)) \cup \partial^+ N^4 (\Gamma (\infty)), \mbox{with} \ \partial^- N^4 (\Gamma (\infty)) \cap \partial^+ N^4 (\Gamma (\infty)) = \Sigma_{\infty}^2 \, .
\end{equation}
We make precise the distinction between $\partial^- N^4$ and $\partial^+ N^4$ by specifying that the $\beta$'s are now 3-balls living, together with the now embedded system $\underset{\ell}{\sum} \, [c_{\ell} (b \, {\rm or} \, r)]$
 hooked at them, entirely inside ${\rm int} \, \partial^+ N^4 (\Gamma (\infty))$, while, for the time being at least, the rest of the link diagram lives entirely inside ${\rm int} \, \partial^- N^4 (\Gamma (\infty))$.
 
 \smallskip
 
 With apropriate framings, this is enough for reconstructing
 $$
 \Theta^4 ((X_0^2 \times r) \cup \Gamma (\infty) \times [r,\beta]) , \varphi) \underset{\rm DIFF}{=} N^4 (X_0^2) \, ,
$$
but we can do better than that too. Starting from the 3-balls $\beta$ inside $\partial^+ N^4 (\Gamma (\infty))$ we get back the whole $N^4 (2\Gamma (\infty))$, which comes now with a splitting by {\it the same} surface $\Sigma_{\infty}^2$,
\begin{equation}
\label{eq4.11}
\partial N^4 (2\Gamma (\infty)) = \partial^- N^4 (2\Gamma (\infty)) \cup \partial^+ N^4 (2\Gamma (\infty)) \, ,
\end{equation}
with $\partial^- N^4 (2\Gamma (\infty)) = \partial^- N^4 (\Gamma (\infty))$, but with a much larger $\partial^+ N^4$.

\smallskip

We have now a grand link, with the $\eta ({\rm green})$ thrown in too, for further purposes, coming with the following {\ibf normal confinement conditions}
\begin{equation}
\label{eq4.12}
\sum_{\ell} c_{\ell} (r \, {\rm or} \, b) + \sum_{\ell} \eta_{\ell} + \sum_1^P \eta_i ({\rm green}) \subset {\rm int} \, \partial^+ N^4 (2\Gamma (\infty)) \, , \ \sum_{i=1}^{\bar P} \Gamma_i + \sum_{j=1}^{\infty} C_j + \sum_{k=1}^{\infty} \gamma_k^0 \subset {\rm int} \, \partial^- N^4 (2\Gamma (\infty)) \, .
\end{equation}
Here, the $c_{\ell} (r)$,  $\Gamma_i$, $C_j$, $\eta_{\ell}$ are attaching zones of internal 2-handles of $N^4 (2X_0^2)$, they have canonical framings, and via all this we can reconstruct $N^4 (2X_0^2)$. The $\eta_i ({\rm green})$ themselves bound exterior discs $\delta_i^2$, and we will not focus here and now on the close connection which $\eta ({\rm green})$ and/or $\delta^2$ may have, or may have had with $c(b)$ (and even with $\gamma^0$). Now, in real life, we will need a certain finite change of the normal confinement conditions (\ref{eq4.12}). There will be a finite system of curves, called generically $\overline{C_f}$, with
$$
\sum \overline{C_f} \subset \sum_1^{\bar P} \Gamma_i + \sum_1^{\infty} C_j \, ,
$$
which are UP at all their corners $P$, and which will be moved isotopically from $\partial^- N^4 (2\Gamma (\infty))$ to $\partial^+ N^4 (2\Gamma (\infty))$. The reasons for this change
\begin{equation}
\label{eq4.13}
(\bar C \subset \partial^- N^4 (2\Gamma (\infty))) \Rightarrow (\bar C \subset \partial^+ N^4 (2\Gamma (\infty)))
\end{equation}
will soon become clear. But the point here is that (\ref{eq4.13}) is actually a {\ibf forced} transformation, the $\overline{C_f}$ may not a priori be UP at all its $P$'s and, anyway, some {\it global} measures will be necessary in order to be able to fit (\ref{eq4.13}) into our whole machinery. In particular, an infinite, locally fine subdivision of $X^2$, {\ibf before doubling}, coming with a {\ibf complete redefinition} of the BLUE labelling, will be needed. What will make this kind of thing possible is the basic commutativity property specific for 2-dimensional collapsing: after any arbitrary collapse, a collapsible 2-dimensional complex stays collapsible. 

\medskip

\noindent IMPORTANT REMARK. As one has already seen, there are quite a number of successive steps in our approach. It is of paramount importance to keep them in correct order, like for instance performing old $\Rightarrow$ new before (\ref{eq4.13}), and (\ref{eq4.13}) itself before doubling$\ldots$

\medskip

The (\ref{eq4.13}) has, of course, to stay compatible with the rest of our construction. Once it is performed, it leads to the real life, forced confinement conditions, which will supersede (\ref{eq4.12}), from now on. They are the following
\begin{equation}
\label{eq4.14}
\sum_{\ell} c_{\ell} (b \, {\rm or} \, r) + \sum_j \eta_j + \sum_f \overline{C_f} \subset \partial^+ N^4 (2\Gamma (\infty)) \supset \sum_1^P \eta_{\alpha} ({\rm green}) \, , 
\end{equation}
$$
\left( \sum_i \Gamma_i + \sum_j C_j - \sum_f \overline{C_f} \right) + \sum \gamma_k^0 \subset \partial^- N^4 (2\Gamma (\infty)) \, .
$$
We will call this from now on the LINK, i.e. the set of internal attaching curves occuring on the LHS of the formula above. There will never be any other violations of the final confinement above, the splitting is sacro-sancted, and nothing is ever allowed to cross $\Sigma_{\infty}^2$.

\smallskip

With this we can open a small prentice, going back to (\ref{eq3.44.1}) which can be explained now. Imagine we would perform a RED diagonalization which would lead to a system of curves
$$
\Gamma_{\bar n + 1} = C(1) \, , \ \Gamma_{\bar n + 2} = C(2) , \ldots , \Gamma_{P+(\bar n - n)} = C(P-n)
$$
which would be attached directly to $\Gamma (3)$. This hypothetical RED diagonalization would have to use both curves $\bar C$ and $(C - \bar C)$, involving thereby serious tresspassing through $\Sigma_{\infty}^2$. This contradicts the sacro-santed principles and, hence, it is {\ibf forbidden}. This proves (\ref{eq3.44.1}).

\smallskip

We go back now to the map
$$
\sum_1^P \delta_i^2 \overset{F}{\longrightarrow} 2X_0^2 \subset \Delta_1^4
$$
from Lemma 9 (see (\ref{eq3.32.1})). This $F$ admits a not everywhere well-defined lift to $\partial N^4 (2X_0^2)$, denoted with the same letter, which we have to look inside, a bit closer, now.

\smallskip

To begin with, there is a piece which we call body $\delta_i^2 \subset \delta_i^2$, and which via the lift of $F$ to $\partial N^4 (2X_0^2)$ goes into $\partial N^4 (2\Gamma (\infty))$. Of course, $\eta_i ({\rm green}) \subset \partial \, {\rm body} \, \delta_i^2$ and also, $\delta^2 - {\rm body} \, \delta^2$ gets  nicely embedded by $F$ into the various lateral surfaces of the 2-handles of $N^4 (2X_0^2)$. We have
\begin{eqnarray}
\label{eq4.15}
&&\partial \, {\rm body} \, \delta_i^2 = \eta_i ({\rm green}) + \sum \{\mbox{the various attaching zones,} \\ 
&&\mbox{call them generically $C_i$, which are such that $F \delta_i^2$ uses $D^2 (C_i)\}$.} \nonumber
\end{eqnarray}
The intersecting part of $F$ is a {\it not everywhere well-defined} immersion, denoted again by the same letter
\begin{equation}
\label{eq4.16}
\sum {\rm body} \, \delta_i^2 \overset{F}{\longrightarrow} \partial N^4 (2\Gamma (\infty)) \, .
\end{equation}
Notice that, the $C_i$ in (\ref{eq4.15}) is part of our LINK, coming with $C_i \subset \partial N^4 (2\Gamma (\infty))$. Careful here, this ``$C_i$'' is just a generic notation for a curve which may be an honest $C \subset \partial^- N^4 (2\Gamma (\infty))$, but which might well be, also, a $\Gamma$ or a $\bar C \subset \partial^+ N^4 (2\Gamma (\infty))$. With all this, the spots where $F$ (\ref{eq4.16}) is {\ibf not} really well-defined correspond to the transversal contacts $F ({\rm body} \, \delta_i^2) \cap C_j$, which we will call {\ibf punctures}. With all these things, here is a typical accident situation, and the description below is supposed to supersede the toy-model considered in connection with (\ref{eq4.2}), in the beginning of this section,
\begin{equation}
\label{eq4.17}
\mbox{Let $L=F ({\rm body} \, \delta_i^2) \cap F ({\rm body} \, \delta_j^2)$ be a {\ibf clasp}, bounded by two punctures}
\end{equation}
$$
p_1 \in F ({\rm body} \, \delta_i^2) \cap C_j \quad {\rm and} \quad p_2 \in F ({\rm body} \, \delta_j^2) \cap C_i \, .
$$
In this situation, we also have, once one goes to dimension four, two transversal contacts
$$
z_1 \in {\mathcal J} \delta_i^2 \cap (D^2 (C_j) \subset 2X_0^2) \, , \ z_2 \in {\mathcal J} \delta_j^2 \cap (D^2 (C_i) \subset 2X_0^2) \, ,
$$
living over $p_1 , p_2$ respectively, as well as a double point
$$
x \in {\mathcal J} \delta_i^2 \cap  {\mathcal J} \delta_j^2 \, .
$$
Figure~2, which lives in a 2-dimensional section through $\Delta_1^4$, may help understand the yoked system of accidents $(z_1 , z_2 , x)$ from (\ref{eq4.17}). Rather than what we had in the context of the toy model (\ref{eq4.2}), the present $p_1 , p_2$ live now at two distinct endpoints of an edge of $2\Gamma (\infty)$, which we call again $p_1 , p_2$. The $F$ (\ref{eq4.16}) is a generic immersion with
$$
M^3 (F) = \phi \ {\rm and} \ FM^2 (F) = \{{\rm clasps} \} \cap \{{\rm ribbons}\}
$$
which, generally speaking, may create a dense web which is highly connected.

\smallskip

We will be now a bit more specific and discuss at some length a typical harder case of (\ref{eq4.17}) where (after a possible permutation of $(1,2)$) we have two pieces $B_i^2 \subset \delta_i^2$, $d_j^2 \subset \delta_j^2$ such that, in the context of (\ref{eq4.17}) we should have
\begin{equation}
\label{eq4.18}
F ({\rm body} \, \delta_i^2) \mid L \subset F ({\rm body} \, B_i^2) \, , \ F ({\rm body} \, \delta_j^2) \mid L \subset F \, d_j^2 \, .
\end{equation}
We know, already, that the transversal contacts ${\mathcal J} \delta^2 \cap \Delta^2 = {\mathcal J} \delta^2 \cap D^2 (\Gamma)$ can only come from the pieces $d^2 \subset \delta^2$ (and moreover this only in the Schoenflies context), which means that with our present specifications, the $\{$extended cocore $z_1\}$ exists always and for sure. We will use it, like in (\ref{eq4.3}). For the sake of the present exposition, let us also assume that $x$ is killed together with $z_1$, leaving us to deal with $z_2$, afterwards, i.e. now. The point is that, in the context (\ref{eq4.18}) we will also have, by construction,
\begin{equation}
\label{eq4.19}
p_2 \in C_i = \bar C_i \subset \partial^+ N^4 (2\Gamma (\infty)) \, .
\end{equation}

$$
\includegraphics[width=8cm]{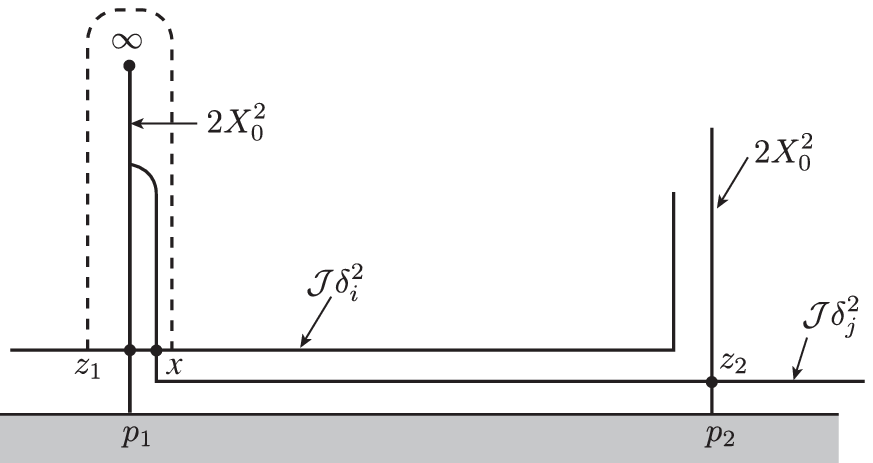}
$$
\begin{figure}[h]
\caption{This is suggesting a 2-dimensional section through the ambient $\Delta_1^4$ showing how the clasp (\ref{eq4.17}) lifts to the yoked system of accidents $(z_1 , x , z_2)$ when we go to four dimensions. The lower hatched zone is supposed to suggest here the $N^4 (2\Gamma (\infty))$, with $L = [p_1 , p_2]$. The arc $[z_1 , \infty]$ is supposed to suggest the compactified $\{$extended cocore $z_1 \}^{\wedge}$. The dotted lines suggest the push of ${\mathcal J} \delta_i^2$ over this $\{$extended cocore $z_1\}^{\wedge}$; it kills here $z_1 + x$. But, in real life, we may well be forced to kill $x$ together with $z_2$, via a completely different procedure which will be outlined in the main text. Our present $z_2$ may or may not possess an $\{$extended cocore $z_2 \}$; in real life at least one of the $z_1 , z_2$ will possess such an extended cocore, anyway. The $\{$extended cocore $z_1 \}^{\wedge} \approx [z_1 , \infty]$ mentioned above, starts inside the $\partial D^2 (C_j)$, which comes with $p_1 \in C_j = \partial D^2 (C_j)$.}
\label{fig1}
\end{figure}

\bigskip

\noindent {\bf Remarks.} A) The reason for (\ref{eq4.13}) was, in retrospect, exactly, to have this (\ref{eq4.19}).

\smallskip

\noindent B) The $d^2$'s only concern $(X_0^2 \times r) \cup (\Gamma (\infty) \times [r,\beta])$ reason why, this piece on which we have focused in (\ref{eq3.26.2}) is more complicated to deal with, than $(\Gamma (\infty) \times [\beta , b]) \cup X_b^2$.

\medskip

Here is the general idea of how one deals with the $z_2$, which lives over $p_2$. We will need an arc $\lambda \subset F \delta_j^2$, joining $p_2$ to some point $q_2 \in \eta_j ({\rm green})$. Eventually, we will want $\lambda$ to live inside $\partial N^4 (2\Gamma (\infty))$, in fact inside $\partial^+ N^4 (2\Gamma (\infty))$, but let use choose to ignore these issues right now, for the purpose of the exposition. The general idea is to use a sliding move of $\eta_j ({\rm green})$ along $\lambda$, until it gets on the other side of $p_2 \in \overline{C_i}$, dragging ${\mathcal J} \delta_j^2$ with it in the process, and thereby destroy the contact $z_2$. In a more precise language, what we mean here is this. Start with the mapping cylinder
$$
{\rm Map} \, ({\mathcal J} \delta_j^2 \approx \delta_j^2 \overset{F}{\longrightarrow} F \delta_j^2 \subset 2X_0^2) \subset \Delta_1^4 \, ,
$$
which is, topologically speaking, essentially $F \delta_j^2 \times [0,\varepsilon]$; then consider a very thin neighbourhood $\lambda \subset U \subset F \delta_j^2$, biting a small arc centered at $q_2$ from $\eta_j ({\rm green})$. Finally, delete the piece of the mapping cylinder living over $U$; this changes $({\mathcal J} \delta_j^2 , F \delta_j^2 , \eta_j ({\rm green}))$ so that $z_2$ disappears.

\smallskip

With this general idea in mind, we go back now to the arc $\lambda$.
\begin{eqnarray}
\label{eq4.20}
&&\mbox{There will be two successive pieces of $\lambda$, first a {\ibf green arc} $\lambda_1 \subset F d_j^2$, connecting $p_2$}  \\
&&\mbox{to some $r_2 \in \gamma_j^0$ and next, a {\ibf dual arc} $\lambda_2 \subset F (\delta_j^2 - d_j^2)$ connecting $r_2$ to $q_2$.} \nonumber \\
&&\mbox{So $\lambda$ takes the form of a paths composition $\lambda = \lambda_1 \cdot \lambda_2$.} \nonumber
\end{eqnarray}

By adding some extra folds to the ``immersion with folds'' $F$, we can arrange this set up in (\ref{eq4.20}), so that
$$
\lambda_1 \subset Fd^2 \subset X_0^2 \times r \, , \ \lambda_2 \subset F (\delta^2 - d^2) \subset \Gamma (\infty) \times [r,b] \cup X_b^2 \, .
$$

The two endpoints $p_2 , q_2$ of $\lambda$ live certainly inside $\partial^+ N^4 (2\Gamma (\infty))$ but, a priori, we may well find that, unless we do something special about it, we have
\begin{equation}
\label{eq4.21}
\lambda \cap \{ \partial^- N^4 (2\Gamma (\infty)) \cup [\mbox{lateral surfaces of the $2$-handles}]\} \ne \emptyset \, .
\end{equation}
We will come back to this unpleasant problem (\ref{eq4.21}) later on. One should remember, at this point, that ${\mathcal J} d^2$, $Fd^2$ were constructed using the 3-dimensional RED collapsing flow, which was still with us up to the $X^2 ({\rm new})$ from (\ref{eq3.18}), but which we have lost by doubling. Now, although this 3-flow has, physically speaking disappeared, its surviving trace on $Fd_j^2$ will be used in order to construct the green arc $\lambda_1$. Also, dually so to say, the 2-dimensional BLUE flow is used for constructing $\lambda_2$. From the begining, the 2-dimensional and the 3-dimensional RED collapsing flows were supposed to be compatible (and we do not explain that now in more detail, the word should suffice here), by doubling we have gained (\ref{eq3.22}), and then, finally the extended cocores use the 2-dimensional RED flow. The net result of all these facts put together, is the following basic item
\begin{equation}
\label{eq4.22}
\{\mbox{extended cocore} \, z_1\}^{\wedge} \cap \lambda = \emptyset \, .
\end{equation}
This means that our two procedures, via which we want to deal with the two ends of the clasp $L$, do not clash with each other. In the same vein, let us notice that, with the green arc $\lambda_1$ guided by the RED flow, and with the dual arc $\lambda_2$ likewise guided by the BLUE flow, if we had not made sure of (\ref{eq3.22}), via the doubling $X^2 \Rightarrow 2X^2$, then in lieu of the normal $\lambda_1 \cap \lambda_2 = \{ r_2 \}$, we would have found, also, plenty of transversal intersections $\lambda_1 \cap \lambda_2$, with disastrous results. As a more general comment, the doubling process (and actually the whole sequence (\ref{eq3.19})) seems to be an essential ingredient for dealing with the accidents.

\smallskip

A lot of various hot issues concerning the accidents have hardly been mentioned so far, and the only thing we can do now is to list at least some of them.
\begin{equation}
\label{eq4.23.1}
\mbox{We certainly want to get rid of (\ref{eq4.21}), with which we cannot live, and achieve}
\end{equation}
$$
\lambda \subset \partial^+ N^4 (2\Gamma (\infty)) \, ,
$$
instead. Here there will be two distinct procedures, one for $\lambda_1$ and another one for $\lambda_2$. For $\lambda_2$, the only issue is to avoid the lateral surfaces of 2-handles. This is an easier issue which can be dealt with by some apropriate subdivisions performed, this time, after doubling, at the level of $X_b^2$ above. [Any pre-doubling subdivision gets automatically ``doubled'' too.] The issue of moving $\lambda_1$ into $\partial^+ N^4 (2\Gamma (\infty))$ is considerably harder and requires some acrobatics which we do not explain here.
\begin{eqnarray}
\label{eq4.23.2}
&&\mbox{One of the offshots of the compatibility between the RED 2-dimensional and 3-dimensional} \\
&&\mbox{flows, will be that $\lambda_1 \cap h = \emptyset$. But we will have contacts $\lambda_1 \cap B \ne \emptyset$. These may threaten} \nonumber \\
&&\mbox{the highly sacro-sancted condition $\eta \cdot B = {\rm id} + {\rm nil}$, and so they are dangerous. They need} \nonumber \\
&&\mbox{hence a treatment, which will not be explained here.} \nonumber
\end{eqnarray}
After all these things, there is no special issue concerning $\lambda_2 \cap (R \cup B)$.
\begin{eqnarray}
\label{eq4.23.3}
&&\mbox{Then, there is also a {\ibf ribbon} analogue of the {\ibf clasp}-accident (\ref{eq4.17}) and this is certainly} \\
&&\mbox{not a trivial thing, contrary to what one may think.} \nonumber
\end{eqnarray}
The correct viewpoint here is to consider the following dense highly connected system (see (\ref{eq4.16}))
$$
FM^2 (F) = \{{\rm clasps}\} \cup \{{\rm ribbons}\} \subset \bigcup_i F ({\rm body} \, \delta_i^2)
$$
which, among other things, raises the hot issue of the unavoidable contacts
$$
\{\mbox{green arcs}\} \cap \{\mbox{clasps and {\ibf ribbons}}\} \ne \emptyset \, ,
$$
into which we will not go here.

\smallskip

It so happens that the big complications of the accidents, are all concentrated along $(X_0^2 \times r) \cup (\Gamma (\infty) \times [r,\beta])$.

\smallskip

We give here a complete description of how accidents ever reach into the region $(\Gamma (\infty) \times [\beta , b]) \cup X_b^2$. This is the following precise local model.
\begin{eqnarray}
\label{eq4.24}
&&\mbox{On the same lines as in (\ref{eq4.17}), we have a clasp $L$, going now along some edge} \\
&&\mbox{$P \times [r,b]$ and involving two $FB^2$'s.} \nonumber
\end{eqnarray}
We have an $FB_i^2$ and an $FB_j^2$, with $z_1$ localized at $P \times r$ and $z_2$ localized at $P \times b$. We treat $(z_1 , x)$ as a single bloc, just like we have done it for (\ref{eq4.17}), except that this is now in earnest, not just an expository pretence. For $z_2$ we use a dual arc $\lambda_2$ confined inside $X_b^2$, without any green arc $\lambda_1$ being necessary here. One can set up things so that there are no $x \in {\mathcal J} M^2 ({\mathcal J})$ localized at $X_b^2$, and all this is more like a simple toy-model of the more difficult case discussed earlier.

\bigskip

\noindent {\bf Remark.} A) It would look, a priori, that when we are dealing with something like (\ref{eq4.17}), we are free to treat $x$ together with $z_1$ {\ibf or} with $z_2$. This is not quite so in real life.

\smallskip

\noindent B) In term of (\ref{eq4.12}), as it stands (and the present discussion is insensitive to the change (\ref{eq4.13})), we consider $({\rm LINK}) \cap \partial^- N^4 (2\Gamma (\infty))$ and its corresponding part of the link projection and link diagram, the only ones which will be discussed now. We know, also, from (\ref{eq4.9}) that
$$
\{\mbox{link diagram}\} = \sum_P \{\mbox{link diagram}\} \mid P \, .
$$
With all this comes now another sacro-sancted principle, which our constructions have always to abide to, namely the following 
\begin{eqnarray}
\label{eq4.25}
&&\mbox{For any individual vertex $P \in \Gamma (\infty)$, inside the corresponding $\{\mbox{link diagram} \, P\} \mid P$, there is} \\
&&\mbox{never an individual line which has both crossings where it is UP and crossings where it is DOWN.} \nonumber
\end{eqnarray}

\noindent C) (A short discussion of (\ref{eq4.25})). So, with (\ref{eq4.25}), any individual line in $\{\mbox{link diagram}\} \mid P$ carries an unmistakable label UP, DOWN, or neutral. Before any (\ref{eq3.19}) is in effect, here is how this could (and will actually) be implemented, at the level of Stage~I, in the previous section. Remembering that $X^2$ is the 2-skeleton of (some cell-decomposition of) $X^4 = X^3 \times R$ and that $X^3$ itself comes with a submersion into $R^3$ (for which, in the context $\Delta^3 \times I$ we have to invoke Smale-Hirsch), we may always assume that, locally at least, the cell-decomposition of $X^4$ is of the form
\begin{equation}
\label{eq4.26}
\{\mbox{a cubically-crystalline decomposition of} \ X^3 \} \times \{\mbox{any subdivision of} \ R\} \, .
\end{equation}
A lot of combinatorial work is required in order to have both $\{$the desirable BLUE and RED features$\}$ AND (\ref{eq4.26}), lumped together inside a single cell-decomposition. But the point here is that with a cell-decomposition like (\ref{eq4.26}), it is not hard to see that (\ref{eq4.25}) is more or less automatically fulfilled. Now, all this was {\it before} we go to the move (\ref{eq3.18}) in Stage~IV (and our (\ref{eq4.25}) which has concerned $[({\rm LINK})$ (\ref{eq4.12})$] \cap \partial^- N^4 (2\Gamma (\infty))$, is insensitive to  whatever may happen in Stage~IV strictly {\it after} the transformation $X^2 ({\rm old}) \Rightarrow X^2 ({\rm new})$). Now, when we go to the real life situation, this passage $X^2 ({\rm old}) \Rightarrow X^2 ({\rm new})$ turns out to be much more complex than what formula (\ref{eq3.18}), as such, may suggest, particularly because {\it we have to abide to} (\ref{eq4.25}). We actually have to use two distinct procedures, once we look into the seams of (\ref{eq3.18}), one for $\Delta^3 \times I$ and then another one for $\Delta^4$ Schoenflies. The difference comes, again, from the existence of the compact product structure, in the first of the two cases.

\smallskip

This is about as much as we will say here, concerning the implementation of (\ref{eq4.25}). We will rather say a few words now about what (\ref{eq4.25}) brings to us. When we consider any $B_i^2$ (\ref{eq3.30}) and we also focus on some $P \in \Gamma (\infty) \times r$ which $B_i^2$ may touch, then $B_i^2 \mid P$ is completely identified by one of the arcs
$$
A \subset \{\mbox{link diagram}\} \mid P \, ,
$$
the connection being that $B_i^2 \mid P = \{$a little triangle spanned by $A$ and the vertex $P \times \beta \}$. When we deal with the accidents, the various $B_i^2 \mid P$ will be (most of the time) dealt with as independent units and, with (\ref{eq4.25}) being satisfied (and also lumping for simplicity's purpose, here, neutral with UP, let us say), each $B_i^2 \mid P$ is
\begin{equation}
\label{eq4.27}
(B_i^2 \mid P) ({\rm UP}) \qquad {\rm OR} \qquad (B_i^2 \mid P)({\rm DOWN}) \, ,
\end{equation}
and never both, simultaneously. Each of these two cases will have to receive a different treatment, reason why we want to keep them distinct. As an illustration for these different treatments, in (\ref{eq4.18}) the $B_i^2 = B_i^2 ({\rm DOWN})$, while at $P \times r$ in (\ref{eq4.24}) we have $B_i^2 ({\rm DOWN}), B_j^2 ({\rm UP})$.

\medskip

\noindent D) Our handling of accidents obviously has to change the topology of the subset
$$
2X_0^2 \cup \sum_1^P {\mathcal J} \delta_i^2 \subset \Delta_1^4 \, ,
$$
but then, in some cases it has to involve changes of the topology of the ambient space $\Delta^4_1$, itself. More explicitely, we may have to use moves which, without touching to $\Delta^4$, of course, locally at least are a brutal change in topology which will turn out, afterwards and this time for global reasons, to leave intact up to diffeomorphism, the pair $(\Delta_1^4 , \Delta^4)$. Without going right now into any particulars, here is how such a {\ibf brutal} move may look like. Consider, in terms of the link diagram, a crossing of two curves none of which are of type $\Gamma_i$. Then interchange the UP/DOWN values at the crossing. Without any loss of generality this does not change the geometric intersection matrices, nor $N^4 (\Delta^2)$ of course. Something quite horrible has, quite clearly, happened locally, but up to diffeomorphism the global topology of the pairs of type
$$
(\Delta_1^4 = N^4 (\Delta^2) \cup \{{\rm collar}\} , N^4 (\Delta^2))
$$
stays intact. So, provided that we do not otherwise conflict with the other structures and/or principles, this is an acceptable move.

\smallskip

With this we close, at the level of the present account, the discussion of the accidents, which we assume, from now on, to have been dealt with already. The last item on our agenda is to give now a very sketchy outline of the grand blue diagonalization. We will be starting now with the diagram (\ref{eq3.51}), and the little BLUE diagonalization (\ref{eq3.52}) is already, and will also constantly be too, with us. Call this the {\ibf initial} level. We have here $\Gamma (3) \subset 2\Gamma (\infty)$, with two families of 1-handles
\begin{equation}
\label{eq4.28}
R({\rm initial}) \subset 2\Gamma (\infty) \supset B \, ,
\end{equation}
where we discard the $R_1 + R_2 + \cdots + R_n$ from (\ref{eq3.3}), as well as the promoted $y_1 + \cdots + y_{P-n} = R_{n+1} + \cdots + R_P$, decreeing that
$$
\sum_1^P b_i \subset R({\rm initial}) \cap B \, ,
$$
which is perfectly legitimate since $\Gamma (3) - \overset{P}{\underset{1}{\sum}} \, b_i$ is now a tree. We do not write $B({\rm initial})$ in (\ref{eq4.28}) since, contrary to what will happen with the $R({\rm initial})$, in the colour-changing process let us call it initial $\Rightarrow$ final, following next, the $B$ will not change at all. We use again the notation
$$
h({\rm initial}) = R({\rm initial}) - \sum_1^P b_i
$$
and, at our present initial stage, we have disjoined partitions
\begin{equation}
\label{eq4.29.initial}
h = (h-B) + h \cap B \, , \ B = (B-h) + B \cap R \, .
\end{equation}
Here $h,R$ are, of course, $h({\rm initial})$, $R({\rm initial})$. With all this, we have $h({\rm initial}) \subset {\rm LAVA} ({\rm initial})$ and we may rewrite (\ref{eq3.48}) as follows
\begin{equation}
\label{eq4.30}
\bar N^4 (\Gamma (3)) = (N^4 (2\Gamma (\infty) - h({\rm initial})) \cup {\rm LAVA}^{\wedge} ({\rm initial}) \, ,
\end{equation}
the two pieces in the RHS being glued along $\delta \, {\rm LAVA} ({\rm initial})$.

\smallskip

In the formula above we apply the prescriptions from Stage~III meaning that
\begin{equation}
\label{eq4.31}
( {\rm LAVA} ({\rm initial}) ,  \delta \, {\rm LAVA} ({\rm initial})) = \left( \bigcup_i  h_i ({\rm initial}) \cup D^2 (C_i) \, , \  \partial \, {\rm LAVA} ({\rm initial}) \cap \partial (N^4 (2\Gamma (\infty) -  h ({\rm initial}))\right) \, .
\end{equation}
The $C \cdot h ({\rm initial})$ is of the easy id $+$ nil form, and hence the pair (\ref{eq4.31}) has the product property. So, at our initial stage we start from
\begin{equation}
\label{eq4.32}
[(N^4 (2\Gamma (\infty)-h ({\rm initial})) \cup {\rm LAVA}^{\wedge} ({\rm initial})] + \sum_1^{\bar P} D^2 (\Gamma_i) \subset \Delta_1^4 \overset{\mathcal J}{\longleftarrow} \sum_1^P \delta_i^2 \, ,
\end{equation}
where the following things happen

\smallskip

a) inside the $[ \ldots \cup \ldots]$, the two corresponding terms are glued together along $\delta \, {\rm LAVA} (\rm {initial})$,

\smallskip

b) between the two compact spaces from the LHS, there is just a collar,

\smallskip

c) we also have
\begin{equation}
\label{eq4.33}
\partial \sum_1^P \delta_j^2 = \sum_1^P \eta_j ({\rm green}) \subset \partial [N^4 (\ldots) \cup {\rm LAVA}^{\wedge} ({\rm initial})] - \sum_1^{\bar P} \Gamma_i \, ,
\end{equation}

d) and finally, apart from (\ref{eq4.33}), ${\mathcal J} \delta^2$ is disjoined from $\Delta^4 \subset \Delta_1^4$, the ${\mathcal J}$ itself being an embedding into this last space. This last point expresses the fact that the accidents are, by now, killed.

\smallskip

Very much like in (\ref{eq4.14}), slightly re-arranged and also considered now at the present initial level after the accidents have been dealt with, we have the following BIG LINK (initial)
\begin{equation}
\label{eq4.34}
\sum_{\ell} c_{\ell} (b \, {\rm or} \, r) + \sum_j \eta_j + \sum_f \overline{C_f} + \sum_1^P \eta_{\alpha} ({\rm green}) \subset \partial^+ N^4 (2\Gamma (\infty)), 
\end{equation}
$$
\left( \sum_{1}^{\bar P} \Gamma_i + \sum_1^{\infty} C_j - \sum_f \overline{C_f} \right) + \sum_1^{\infty} \gamma_k^0 \subset \partial^- N^4 (2\Gamma (\infty)) \, .
$$

To complete the picture at the initial level, let us add the following two items too. The small blue diagonalization (\ref{eq3.52}) is, and will still constantly be, with us from now on. Finally, the only obstruction which has been left on our way now, is the following finite set, which was already identified in (\ref{eq3.53}), namely
\begin{eqnarray}
\label{eq4.34.1}
&&\mbox{the $h_k \in h ({\rm initial}) - B$ such that $h_k \subset \overset{P}{\underset{1}{\sum}} \ \{\mbox{extended cocore} \, b_j \}^{\wedge}$} \\
&&\mbox{(the $1$-handles of $\Delta^4$) and which, also, are touched by $\overset{P}{\underset{1}{\sum}} \ \eta_i ({\rm green})$.} \nonumber
\end{eqnarray}

With all these things, what comes next is a big geometric transformation, which we call
\begin{equation}
\label{eq4.35}
\mbox{The CHANGE OF COLOUR initial $\Rightarrow$ final,}
\end{equation}
at the final level of which we will find a context analogous (modulo some important changes) to the one of the initial level, but where now the GRAND BLUE DIAGONALIZATION IS IN PLACE. The transformation (\ref{eq4.35}) will leave $(2\Gamma (\infty) , B , \eta ({\rm green}))$, eventually, invariant. By ``eventually'' we mean here that the initial and final levels for these objects will be rigorously the same, but with a lot of drastic transformations occuring between. The (\ref{eq4.35}) also comes with a change
\begin{equation}
\label{eq4.36}
\{ h({\rm initial}) , ({\rm LAVA} ({\rm initial}) , \delta \, {\rm LAVA} ({\rm final}))\} \Rightarrow \{h ({\rm final}) , ({\rm LAVA} ({\rm final}) , \delta \, {\rm LAVA} ({\rm final}))\}
\end{equation}
satisfying the usual condition
\begin{equation}
\label{eq4.37}
\delta \, {\rm LAVA} ({\rm final}) =  \partial \, {\rm LAVA} ({\rm final}) \cap \partial (N^4 (2\Gamma (\infty)) - h({\rm final})) \, .
\end{equation}
It will turn out now, that the change of colour process (\ref{eq4.35}) brings with it a quite serious violation of the RED $C \cdot h = {\rm id} + {\rm nil}$ feature. This means that, in order to retain the product property for the pair
$$
({\rm LAVA} ({\rm final}) , \delta \, {\rm LAVA} ({\rm final})) \eqno (*)
$$
we cannot any longer use the exact prescription from the Stage~III. We had used these prescriptions for (\ref{eq4.31}), but for defining correctly $(*)$, some modifications of the standard prescriptions will be needed. Without going into that, right now, with $[N^4 (2\Gamma (\infty) - h ({\rm initial})] \cup {\rm LAVA}^{\wedge} ({\rm initial})$, replaced now by
\begin{equation}
\label{eq4.38}
[N^4 (2\Gamma (\infty)) - h ({\rm final})] \cup {\rm LAVA}^{\wedge} ({\rm final}) \, ,
\end{equation}
and with $(*)$ which still retains the product property, we have, at the final level, a context just like in (\ref{eq4.32}). In particular, the product property, provides us with a properly embedded system of 3-balls
\begin{equation}
\label{eq4.39}
\sum_1^P \{\mbox{extended cocore} \, b_i\}^{\wedge} ({\rm final}) \, .
\end{equation}

The next lemma should explain the term ``change of colour''.

\bigskip

\noindent {\bf Lemme 13.} {\it For any finite subset $S \subset h({\rm initial}) - B$, we can find a larger, still finite subset
\begin{equation}
\label{eq4.40}
S \subset S_1 \subset h({\rm initial}) - B
\end{equation}
for which there exists an
\begin{equation}
\label{eq4.41}
S_2 \subset B - R({\rm initial}) , \ \mbox{with} \ {\rm card} \, S_2 = {\rm card} \, S_1
\end{equation}
such that the following things should happen.}

\medskip

1) {\it If one defines
\begin{equation}
\label{eq4.42}
R({\rm final}) = \sum_1^P b_i + (h({\rm initial}) - S_1) + S_2 , \mbox{and hence also} \ h({\rm final}) = (h({\rm initial}) - S_1) + S_2 \, ,
\end{equation}
then, just like $2\Gamma (\infty) - B$ and $2\Gamma (\infty) - R({\rm initial})$, the $2\Gamma (\infty) - R({\rm final})$ is again a tree.}

\medskip

2) {\it In a manner which the very sketchy proof below will make explicit (at least up to a certain extent), this comes with the items {\rm (\ref{eq4.36}), (\ref{eq4.37}), (\ref{eq4.38}), (\ref{eq4.39})}, i.e. with the following analogue of {\rm (\ref{eq4.32})}, at the final level
\begin{equation}
\label{eq4.43}
[(N^4 (2\Gamma (\infty) - h ({\rm final})) \cup {\rm LAVA}^{\wedge} ({\rm final})] + \sum_1^{\bar P} D^2 (\Gamma_i) \subset \Delta_1^4 \overset{\mathcal J}{\longleftarrow} \sum_1^P \delta_i^2 \, ,
\end{equation}
satisfying the analogue of} (\ref{eq4.33}).

\medskip

3) (PUNCH LINE) {\it If $S$ is large enough so as to contain the finite set from {\rm (\ref{eq4.34.1})}, then we also have}
\begin{equation}
\label{eq4.44}
\sum_1^P \eta_i ({\rm green}) \cap {\rm LAVA}^{\wedge} ({\rm final}) = \emptyset \, .
\end{equation}

\bigskip

Before we go to a very sketchy and impressionistic outline of proof for this last lemma, let us notice that, modulo everything said so far, the (\ref{eq4.44}) should clinch the proofs of the two Theorems~1 and 2. With this we list now the kind of steps invoked in the proof of Lemma~13. This is just a sketchy outline, of course.

\medskip

I) At the 1-dimensional level, the (\ref{eq4.35}) is a sequence of embedded transformations of $N^4 (2\Gamma (\infty))$ inside $\Delta_1^4$, which are 1-handle slides keeping all the time the {\ibf splitting intact} and changing the positions of 1-handle cocores around. In the end (but {\ibf only} in the end), we find ourselves with exactly the same $N^4 (2\Gamma (\infty))$ as in the beginning, and also with the transformation of pairs
$$
(N^4 (2\Gamma (\infty)) , R({\rm initial})) \Rightarrow (N^4 (2\Gamma (\infty)) , R({\rm final})) \, .
$$
But, in this eventual transformation, the $B$ stays put (although it might have done horrible things at intermediary stages, a leitmotif in this present story). So, it is only $h-B$ which actually changes, when we move from ``initial'' to ``final''.

\medskip

II) The 1-handle slides from I), drag along the curves and the 2-handles, internal or external attached along them. The confinement conditions are {\ibf never} violated and, with a lot of intermediary stages we get a transformation at the level of (\ref{eq4.34})
\begin{equation}
\label{eq4.45}
\mbox{BIG LINK (initial) $\Rightarrow$ BIG LINK (final),}
\end{equation}
at the end of which (but only at the end) we find that, as subsets of $2\Gamma (\infty)$, we have the equality
\begin{equation}
\label{eq4.46}
\{\mbox{BIG LINK (initial)}\} \cap \partial^+ N^4 (2\Gamma (\infty)) = \{ \mbox{BIG LINK (final)}\} \cap \partial^+ N^4 (2\Gamma (\infty)) \, .
\end{equation}

III) Our set-up is such that, once we choose to forget the intermediary stages and only look at the initial and final levels, then both the geometric intersection matrices $\eta \cdot B$ and $\eta ({\rm green}) \cdot B$ stay put. But not so on the RED side when we will actually have 
\begin{equation}
\label{eq4.47}
C({\rm final}) \cdot h ({\rm final}) \ne {\rm id} + {\rm nilpotent} \, .
\end{equation}

\medskip

\noindent {\bf Remarks.} Both $C ({\rm initial})$ and $C({\rm final})$ contain curves $\bar C \subset \partial^+ N^4 (2\Gamma (\infty))$. Also, in a sense which is not too hard to make precise the $C ({\rm initial}) \cdot h({\rm initial})$ and $C ({\rm final}) \cdot h({\rm final})$ only differ by a {\ibf finite} matrix, let us say that the {\ibf violation} of the RED feature id $+$ nil displayed above, is compact.

\medskip

IV) Once we have lost the RED id $+$ nil, we can no longer proceed exactly like in the context of Stage~III. So now we will have, by definition
\begin{equation}
\label{eq4.48}
{\rm LAVA} ({\rm final}) = \left[ \sum_i (h_i ({\rm final})) \cup D^2 (C_i ({\rm final})) \right] \cup \{\mbox{some {\ibf additional} pieces which we call {\ibf lava bridges}}\}.
\end{equation}

The $\sum$ (lava bridges) is compact and its role is to make that the following should happens
\begin{eqnarray}
\label{eq4.49}
&&\mbox{The pair $({\rm LAVA} ({\rm final}) , \delta \, {\rm LAVA} ({\rm final}))$ continues to have the product property,} \\
&&\mbox{actually with the same lamination as before, i.e. ${\mathcal L} ({\rm initial}) = {\mathcal L} ({\rm final})$.} \nonumber
\end{eqnarray}

V) There is a geometric transformation
\begin{equation}
\label{eq4.50}
({\rm LAVA} ({\rm initial}) , \delta \, {\rm LAVA} ({\rm initial})) \Rightarrow ({\rm LAVA} ({\rm final}) , \delta \, {\rm LAVA} ({\rm final}))
\end{equation}
which has the virtue that its various intermediary stages exhibit, explicitely, the conservation of the product property. Bot the initial and final levels of (\ref{eq4.50}) are naturally embedded inside the ambient $\Delta_1^4$, but {\ibf not} so all the intermediary steps. Let us say here that (\ref{eq4.50}) is an ``{\ibf allowable lava move}''. Also, in order for the product property to be preserved, notwithstanding the fact that we only have a compact violation of id $+$ nil in the context of (\ref{eq4.50}), some global conditions, involving the totality of $N^4 (2X_0^2)$ will have to be paid attention too.

\medskip

VI) We will constantly have
\begin{equation}
\label{eq4.51}
\{\mbox{lava bridges}\} \cap \partial N^4 (2\Gamma (\infty)) \subset \partial^- N^4 (2\Gamma (\infty)) \, ,
\end{equation}
far from $\eta ({\rm green})$.

\medskip

VII) Here is a heuristic argument suggesting why we should have (\ref{eq4.44}) (the punch line). We already know that
$$
\eta ({\rm green}) \cdot \left( B - \sum_1^P b_j \right) = 0
$$
and, quite clearly the $S_2$ (\ref{eq4.41}) $\subset \, B - \overset{P}{\underset{1}{\sum}} \, b_i$. Assume now that $S \supset \{\mbox{the finite set of $h_k$'s from (\ref{eq4.34.1}), call it}$ $S_3\}$. It is also, only through $S_3$ that $\eta ({\rm green})$ touches LAVA (initial). So, if $S_3$ is changed into a piece of the BLUE $S_2$, then we get (\ref{eq4.44}); end of proof! The real life argument is, of course, a bit more complex, but this is, anyway, the idea.

\medskip

VIII) All the steps outlined so far may seem a bit mysterious, and so I would like to focus now, for a minute or so, on the exact moment when the actual change of colour takes place, and in particular on the geometry which comes with it. We are supposed to be now somewhere in the middle of the process (\ref{eq4.35}), at a time which I will call ``intermediary''. This comes with a $2\Gamma (\infty)$ (intermediary) quite different from $2\Gamma (\infty)$, but which is split along an $\Sigma_{\infty}^2$ (intermediary) as follows
\begin{equation}
\label{eq4.52}
\partial N^4 (2\Gamma (\infty) ({\rm intermediary})) = \partial^+ N^4 \cup \partial^- N^4 , \partial^- N^4 \cap \partial^+ N^4 = \Sigma_{\infty}^2 ({\rm intermediary}) \, .
\end{equation}
There is also a BIG LINK (intermediary) satisfying the obvious confinement condition and an $R$ (intermediary), containing some precise, interesting RED element
\begin{equation}
\label{eq4.53}
h_i \in (R({\rm intermediary}) - B) \cap S_1
\end{equation}
which is such that the time has come for trading it for its BLUE counterpart $B_i \in S_2 - R$ (intermediary). And it is the geometry of this trading step $h_i \leftrightarrow B_i$, which we want to explain now. But before we can do this we have to start by unravelling one of the main virtues of the step (I) above. Some notations will be necessary here. Let us consider the 3-ball $B^3$, together with the splitting of its boundary by the equatorial circle, call it
$$
\partial B^3 = \partial^+ B^3 \cup \partial^- B^3
$$
where $\partial^{\pm} B^3$ are the two hemispheres. Next, consider a long cylinder $B^3 \times [-N,N]$, with $[-N,N] \subset \{\mbox{some $x$-axis}\}$, and along this $x$-axis we consider the four quantities
$$
-N < r < b < N \, .
$$
Notice the following splitting for the lateral surface of our cylinder
\begin{equation}
\label{eq4.54}
\partial B^3 \times [-N,N] = (\partial^- B^3 \times [-N,N]) \cup (\partial^+ B^3 \times [-N,N]) \, .
\end{equation}
With all these things, what step (I) does for us, is to generate an embedding
\begin{equation}
\label{eq4.55}
(B^3 \times [-N,N] , \partial B^3 \times [-N,N] ) \overset{\ell}{\longrightarrow} (N^4 (2\Gamma (\infty) ({\rm intermediary})) , \partial N^4 (2\Gamma (\infty) ({\rm intermediary})) \, ,
\end{equation}
which is such that
\begin{equation}
\label{eq4.55.1}
\ell (B^3 \times r) = h_i \, , \ \ell (B^3 \times b) = B_i \ \mbox{and, moreover} \ {\rm Im} \ell \cap (B \cup R) = \{ h_i , B_i \} + \{\mbox{some harmless} \, B-R\} \, .
\end{equation}
\begin{equation}
\label{eq4.55.2}
\mbox{The embedding $\ell$ is compatible with the splittings (\ref{eq4.54}) (at the source) and (\ref{eq4.52}) (at the target), i.e.}
\end{equation}
$$
\ell (\partial^{\pm} B^3 \times [-N,N]) \subset \partial^{\pm} N^4 (2\Gamma (\infty) ({\rm intermediary})) \, .
$$
We can omit to write the ``$\ell$'' explicitely, from now on. At the intermediary moment where we find ourselves now, the (\ref{eq4.46}) is, most likely, violated. What we have, instead, are the following two sets, to consider next
\begin{equation}
\label{eq4.56}
\Lambda^+ = \{\mbox{BIG LINK (intermediary)}\} \cap \partial^+ B^3 \times [-N,N] \ {\rm and}
\end{equation}
$$
\Lambda^- = \{\mbox{BIG LINK (intermediary)} + [\mbox{lava bridges}]\} \cap \partial^- B^3 \times [-N,N] \, .
$$
Here are some comments concerning (\ref{eq4.55}) and (\ref{eq4.56}).
\begin{eqnarray}
\label{eq4.56.1}
&&\mbox{We have $B^3 \times r \subset {\rm LAVA}$ , $B^3 \times b \not\subset {\rm LAVA}$, the lava under discussion now,} \\
&&\mbox{ being the one at the intermediary moment before any action.} \nonumber
\end{eqnarray}
\begin{equation}
\label{eq4.56.2}
\mbox{Except for $(\eta ({\rm green}) + c(b)) \cap \Lambda^+$ and for $(\Gamma + \gamma^0) \cap \Lambda^-$, everything else in $\Lambda^{\pm}$ is lava.}
\end{equation}
\begin{equation}
\label{eq4.56.3}
\mbox{Consider any $h_u$ and any curve $C_v$ which is lava (which, in terms of (\ref{eq4.34}) may mean (actual) $C,\eta$ or $c(r)$).}
\end{equation}

\noindent Any contact $C_v \cdot h_u$, a priori {\ibf sticks}, in the sense that if we severe it, then we might well destroy the product property of lava. Now, it so happens that our LAVA has a bit more internal structure than what has been displayed, so far. One of the consequences of this not yet explained additional structure, is that, if $C_v \subset \partial^+ N^4 (2\Gamma (\infty))$ and if $h_u$ pertains to $X_0^2 \times r$,  then the contacts $C_v \cdot h_u$ do not stick, one can sever them without destroying the product property.

\smallskip

With all this we consider now an internal transformation $T$ of $N^4 (2\Gamma (\infty) ({\rm intermediary}))$, operating as follows.
\begin{eqnarray}
\label{eq4.57.1}
&&\mbox{The transformation $T$ applied to the space $N^4 (2\Gamma (\infty) ({\rm intermediary}))$ is a simple isotopic} \\
&&\mbox{diffeomorphism, respecting the splitting, and having all of its active part concentrated} \nonumber \\ 
&&\mbox{inside $B^3 \times [-N,N]$.} \nonumber
\end{eqnarray}
\begin{eqnarray}
\label{eq4.57.2}
&&\mbox{If $x$ is the coordinate along $[-N,N]$, then $T$ is the identity on the factor $B^3$, while along} \\
&&\mbox{$[-N,N]$ it is the translation $T(x) = x + (b-r)$, dampened so that it becomes the identity,} \nonumber \\
&&\mbox{again, in the neighbourhood of $\pm N$.} \nonumber
\end{eqnarray}
\begin{equation}
\label{eq4.57.3}
\mbox{So, geometrically speaking} \ T (B^3 \times r) = B^3 \times b \, ,
\end{equation}
which we accompany by the following decrees. To begin with, we change $R ({\rm intermediary})$ into $R$ (intermedia\-ry) $ - \{ h_i \} + \{ B_i \}$ deciding, also that now $B_i \in R \cap B$. Next, we also decree that
$$
B^3 \times b_i \subset LAVA \, ,
$$
and we completely discard the $h_i$, as such, from the rest of our procedure.
\begin{eqnarray}
\label{eq4.57.4}
&&\mbox{Declaring that $B^3 \times b$ is LAVA, is a less innocent operation than it may look at first sight} \\
&&\mbox{since $B^3 \times b$ certainly comes with contacts $(B^3 \times b) \cap \Lambda^{\pm}$ of its own, which we have to worry} \nonumber \\
&&\mbox{about now.} \nonumber
\end{eqnarray}
But let us first describe the action of $T$ on $\Lambda^{\pm}$, which will be done in the (\ref{eq4.57.5}) below.
\begin{eqnarray}
\label{eq4.57.5}
&&\mbox{({\ibf The main step}) On the side of $\partial^- B^3 \times [-N,+N]$, we let $\Lambda^-$ go solidarily with $B^3 \times r$,} \\
&&\mbox{i.e. we apply to it the same geometrical move $x \mapsto x+(b-r)$ like in (\ref{eq4.57.2}) above. But} \nonumber \\
&&\mbox{then on the $\partial^+ B^3 \times [-N,+N]$ side we take $T \mid \Lambda^+ = {\rm identity}$, i.e. {\ibf we leave $\Lambda^+$ in place},} \nonumber \\
&&\mbox{\ibf without budging it.} \nonumber
\end{eqnarray}

All this requires some explaining. The $B^3 \times r$ has taken along with it, to its new position $B^3 \times b$, all the contacts with $\Lambda^-$ which it had. This includes, of course lava, which comes now on top of $B^3 \times b$, i.e. new lava connections. But that is fine, since $B^3 \times r$ has been just changed into $B^3 \times b$. Also, at the same time, the same $x \mapsto x + [b-r]$ removes all the old connections $(B^3 \times b) \cap \Lambda^-$ which might have been there. One can easily see that, at the local level of $\partial^- B^3 \times [-N,N]$, all this is OK, lava-wise. Also, we have dragged along the non lava part of $\Lambda^-$, which inside its confinement site $\partial^- N^4$ is, generally speaking, entangled with the rest of $\Lambda^-$. This way, we have avoided the danger of tearing apart the topology of
$$
[N^4 (2\Gamma (\infty) - h) \cup {\rm LAVA}^{\wedge}] + \sum D^2 (\Gamma) \, .
$$

This is all we have to say on the $\partial^- N^4$ side. On the $\partial^+ N^4$ side, the lava connections coming from $(B^3 \times r) \cap \Lambda^+$ have been severed, but this is allowable, via (\ref{eq4.56.3}). Then also, new lava connections coming with $(B^3 \times b) \cap \Lambda^+$ have been established and this other brutal move is also allowable, because of another global property of the lava, dual to (\ref{eq4.56.3}), and quite similar to it. So, our whole movement $T$ which has exchanged $h_i \in S_1$ with $B_i \in S_2$ has kept the product property of lava intact, while at the same time, achieving the following basic thing 
\begin{equation}
\label{eq4.58}
\mbox{Any contact $\eta ({\rm green}) \cdot h_i$ has gone and no contact $\eta ({\rm green}) \cdot B_i$ has appeared instead.}
\end{equation}

This is obviously the kind of thing we need, for getting (\ref{eq4.44}).

\bigskip

\noindent {\bf A final remark.} Notice that it is the splitting $+$ confinement, both sacro-sancted principles in this paper, which allow us to operate independently on $\Lambda^+$ and $\Lambda^-$ without getting them entangled with each other. But then, splitting $+$ confinement are necessary all over the place in the proof of Lemma~13, like for instance for restoring, at the final level, the
$$
(N^4 (2\Gamma (\infty)) , B , \eta ({\rm green}) (\mbox{together with the (BIG LINK)} \, \cap \partial^+ N^4)) \, ,
$$
exactly as they were at the initial one.


\newpage

\begin{thebibliography}{999}
\bibitem{Ga} {\sc D. Gabai}, Valentin Po\'enaru's Program for the Poincar\'e Conjecture, in the volume {\it Geometry Topology and Physics for Raoul Bott} (ed. by S.T.~Yau), International Press, pp.~139-169 (1994).
\bibitem{Ma} {\sc B. Mazur}, On embedding of spheres, {\it BAMS}, {\bf 65}, pp.~59-65 (1959).
\bibitem{Mo} {\sc J.W. Morgan}, Recent progress on the Poincar\'e Conjecture and the classification of 3-manifolds, {\it BAMS}, {\bf 42}, pp.~57-78 (2004).
\bibitem{PoI} ([PoI]) {\sc V. Po\'enaru}, The collapsible pseudo-spine representation theorem, {\it Topology}, vol.~31, {\bf 3}, pp.~625-636 (1992).
\bibitem{PoII} ([PoII]) {\sc V. Po\'enaru}, Infinite processes and the 3-dimensional Poincar\'e Conjecture, II: The Honeycomb representation theorem, {\it Pr\'epublications d'Orsay}, 93-14 (1993).
\bibitem{PoIII} ([PoIII]) {\sc V. Po\'enaru}, Infinite processes and the 3-dimensional Poincar\'e Conjecture, III: The algorithm, {\it Pr\'epublications d'Orsay}, 92-10 (1992).
\bibitem{PoIV-A} ([PoIV-A]) {\sc V. Po\'enaru}, Processus infini et conjecture de Poincar\'e en dimension trois, IV: Le th\'eor\`eme de non sauvagerie lisse (The smooth tameness theorem), Part A, {\it Pr\'epublications d'Orsay}, 93-83 (1992).
\bibitem{PoIV-B} ([PoIV-B]) {\sc V. Po\'enaru}, Processus infini et conjecture de Poincar\'e en dimension trois, IV: Le th\'eor\`eme de non sauvagerie lisse (The smooth tameness theorem), Part B, {\it Pr\'epublications d'Orsay}, 95-33 (1995).
\bibitem{PoV} This [PoV], {\it Pr\'epublications d'Orsay}, 94-25 (1994) is completely superseded by \cite{PoV-A}, \cite{PoV-B}.
\bibitem{PoV-A} ([PoV-A]) {\sc V. Po\'enaru}, Geometric simple connectivity in four-dimensional differential topology, PartA, {\it IHES Pr\'epublications} M/01/45 (2001), http://www.ihes.fr/PREPRINTS.M01/Resu/resu-M01-45.html
\bibitem{PoV-Aoutline} {\sc V. Po\'enaru}, Geometric simple connectivity in four-dimensional differential topology: An outline, Preprint Trento Univ. UTM 649 (2003), http://eprints.biblio.unitn.it/archive/00000660/02/UTM649-pdf
\bibitem{PoV-B} ([PoV-B]) {\sc V. Po\'enaru}, Manuscript.
\bibitem{PoVI} ([PoVI]) {\sc V. Po\'enaru}, The strange compactification theorem, Part A, {\it IHES Pr\'epublications} M/95/15 (1995); Part B, {\it IHES Pr\'epublications} M/96/43 (1996); Part C, {\it IHES Pr\'epublications} M/97/43 (1997); Part D, {\it IHES Pr\'epublications} M/97/59 (1997); Part E is in process of being typed at IHES.
\bibitem{Po-B} {\sc V. Po\'enaru}, A program for the Poincar\'e Conjecture and some of its ramifications, in the volume {\it Topics in Low-Dimensional Topology} (ed. A.~Banyaga, H.~Movahedi-Lankarani, R.~Wells), World Scientific, pp.~65-88 (1999).
\bibitem{Po-S} {\sc V. Po\'enaru}, Geometric Simple Connectivity and Low-Dimensional Topology, {\it Proceedings of the Steklov Institute}, {\bf 247}, pp.~195-208 (2004).
\bibitem{Po-T} {\sc V. Po\'enaru}, Three lectures on higher-dimensional methods in three-dimensional topology, {\it Proceedings of the F.~Tricerri Memorial Conference}, Suppl. ai Rendiconti del Circolo Matematico di Palermo, S-II N49, pp.~203-217 (1997).
\bibitem{W1} {\sc J.H.C. Whitehead}, On $C^1$-complexes, {\it Ann. of Math.} {\bf 41}, pp.~809-824 (1940).
\bibitem{W2} {\sc J.H.C. Whitehead},  A certain open manifold whose group is unity, {\it Q. J. of Math.}, {\bf 6}, pp.~268-279 (1935).
\end{thebibliography}
\end{document}